\documentclass[10pt,twoside, a4paper, english, reqno]{amsart}
\usepackage[dvips]{epsfig}
\usepackage{amscd}
\usepackage{a4wide}
\usepackage{amssymb}
\usepackage{amsthm}
\usepackage{amsmath}
\usepackage{bm} 
\usepackage{latexsym}
\usepackage{enumitem}
\usepackage{tabularx}
\usepackage{booktabs} 
\usepackage{siunitx}
\usepackage{array} 
\numberwithin{figure}{section}
\numberwithin{table}{section}
\usepackage{graphicx,caption}
\DeclareMathOperator{\sech}{sech}

\usepackage{calc}
\usepackage{appendix}
\usepackage{graphicx}
\usepackage{amsmath,amsthm,bm,mathrsfs}
\usepackage[utf8]{inputenc} 
\usepackage{url}            
\usepackage{hyperref}

\usepackage[foot]{amsaddr}

\usepackage{caption}    

\usepackage{xcolor}
\definecolor{mygreen}{rgb}{1,0.3,0.6}
\definecolor{mypurple}{HTML}{800080}
\definecolor{myorange}{HTML}{FFA500}
\definecolor{myteal}{HTML}{008080}
\definecolor{mybrown}{HTML}{8B4513}

\hypersetup{
    colorlinks=true,
    linkcolor=blue,
    filecolor=red,      
    urlcolor=blue,
    citecolor=red,
    }
\usepackage{amsfonts}       
\usepackage{nicefrac}       
\usepackage{microtype}      
\theoremstyle{plain}
\usepackage{doi}

\usepackage{upref}
\usepackage{hyperref}
\usepackage{ulem}
\usepackage{enumitem}
\usepackage{graphicx}
\newcommand{\norm}[1]{\left\Vert#1\right\Vert}

\newcommand{\R}{\mathbb R}
\newcommand{\Z}{\mathbb{Z}}
\newcommand{\N}{\mathbb N}

\newtheorem{theorem}{Theorem}[section]
\newtheorem{proposition}[theorem]{Proposition}
\newtheorem{lemma}[theorem]{Lemma}

\newtheorem{remark}[theorem]{Remark}

\numberwithin{equation}{section}     
\numberwithin{figure}{section}
\numberwithin{table}{section}

\allowdisplaybreaks

\newcounter{asnr}

{\ifnum\value{asnr}=0 \stepcounter{asnr} 
  \begin{enumerate}[label=\textbf{A}.\arabic{enumi}]
    \else
    \begin{enumerate}[label=\textbf{A}.\arabic{enumi},resume] \fi}
{\end{enumerate}}

\newcounter{defnr}

{\ifnum\value{defnr}=0 \stepcounter{defnr} 
  \begin{enumerate}[label=\textbf{D}.\arabic{enumi}]
    \else
    \begin{enumerate}[label=\textbf{D}.\arabic{enumi},resume] \fi}
{\end{enumerate}}

\numberwithin{equation}{section} \allowdisplaybreaks

\title[Spectral method for zero dispersion limit]
{Spectral Galerkin method for the zero dispersion limit of the fractional Korteweg-de Vries equation}
\date{}

\author{\textsuperscript{a}Mukul Dwivedi\textsuperscript{*} \and \textsuperscript{b}Tanmay Sarkar}

\address{\textsuperscript{a,b}Department of Mathematics, Indian Institute of Technology Jammu, 
Jagti, NH-44 Bypass Road, Post Office Nagrota, Jammu - 181221, India.}

\email{\textsuperscript{*}Corresponding author: 2020rma1031@iitjammu.ac.in}



\thanks{}

\begin{document}

\begin{abstract}
We present a fully discrete Crank-Nicolson Fourier-spectral-Galerkin (FSG) scheme for approximating solutions of the fractional Korteweg-de Vries (KdV) equation, which involves a fractional Laplacian with exponent $\alpha \in [1,2]$ and a small dispersion coefficient of order $\varepsilon^2,~\varepsilon\ll 1$. The solution in the limit as $\varepsilon \to 0$ is known as the zero dispersion limit. We demonstrate that the semi-discrete FSG scheme conserves the first three integral invariants, thereby structure preserving, and the fully discrete FSG scheme is $L^2$-conservative, ensuring stability. Using the compactness argument, we constructively prove the convergence of the approximate solution to the unique solution of the fractional KdV equation in $C([0,T]; H_p^{1+\alpha}(I))$ for the periodic initial data in $H_p^{1+\alpha}(I)$. The devised scheme achieves spectral accuracy for the initial data in $H_p^r,~r \geq 1+\alpha$ and exponential accuracy for the analytic initial data.

Additionally, we establish that the approximation of the zero dispersion limit obtained from the fully discrete FSG scheme converges to the solution of the Hopf equation in $L^2$ as $\varepsilon \to 0$, up to the gradient catastrophe time $t_c$. Beyond $t_c$, numerical investigations reveal that the approximation converges to the asymptotic solution, which is weakly described by the Whitham’s averaged equation within the oscillatory zone for $\alpha = 2$. Numerical results are provided to demonstrate the convergence of the scheme and to validate the theoretical findings.
\end{abstract}

\maketitle

\section{Introduction}

In recent years, differential equations involving the fractional Laplacian operator have become significant tools across various scientific, economic, and engineering fields due to their capability to model processes with anomalous diffusion and non-local interactions \cite{di2012hitchhikers,podlubny1998fractional}. The fractional Laplacian, defined as a non-local operator, is particularly effective in capturing localized phenomena, such as image processing, fluid flow in porous media, and complex plasma dynamics (see  \cite{kwasnicki2017ten,podlubny1998fractional} and references therein).

This paper focuses on the Cauchy problem for the fractional Korteweg–de Vries (KdV) equation, defined by
\begin{equation}\label{fkdv}
\begin{cases}
      u_t + 6uu_x - \varepsilon^2 \mathcal{D}^{\alpha} u_x = 0,  & (x,t) \in \R_T := \mathbb{R} \times (0,T], \\
      u(x,0) = u_0(x), & x \in \mathbb{R},
\end{cases}
\end{equation}
where $\mathcal{D}^\alpha := (-\Delta)^{\alpha/2}$ represents the fractional Laplacian with exponent $\alpha \in [1, 2]$, $\varepsilon^2$ is a small dispersion coefficient, $u_0$ is the given initial condition, and $u: \R_T \rightarrow \mathbb{R}$ is the solution to be determined. The fractional Laplacian is defined in terms of the Fourier transform for the exponent $\alpha \in (0, 2]$  as
\begin{equation}\label{fracLF}
    \widehat{[\mathcal{D}^{\alpha} u]}(\xi) = |\xi|^\alpha \widehat{[u]}(\xi),
\end{equation}
where $\widehat{[u]}$ denotes the Fourier transform of $u$. This operator introduces non-local effects into the dynamics of the equation, thereby extending the classical KdV model to incorporate fractional dispersion.

The fractional KdV equation \eqref{fkdv} is a nonlinear, non-local dispersive model that arises in the study of weakly nonlinear long internal waves \cite{dutta2021operator,dwivedi2024fully,kenig1991well, kenig1993cauchy,molinet2018well}. When $\alpha = 2$, the equation reduces to the classical KdV equation, known for modeling solitons and nonlinear wave phenomena \cite{sjoberg1970korteweg,Bona1975KdV, dutta2021operator, korteweg1895xli}. For $\alpha = 1$, it corresponds to the Benjamin-Ono (BO) equation, which describes one-dimensional internal waves in deep stratified fluids \cite{fokas1981hierarchy, kenig1994generalized,tao2004global,molinet2018well}. The well-posedness of the fractional KdV equation has been extensively studied in the Sobolev spaces $H^s(\R)$, with significant results including local well-posedness for $s > 3/4$ (see \cite{kenig1991well, kenig1993cauchy}), and global well-posedness for $L^2$ initial data using frequency-dependent renormalization techniques \cite{herr2010differential}. 

Despite these theoretical advancements, numerical methods for solving the fractional KdV equation \eqref{fkdv} for $\alpha \in (1,2)$ are limited. Existing approaches include operator splitting \cite{dutta2021operator}, Galerkin methods \cite{dwivedi2023stability}, and finite difference schemes \cite{dwivedi2024fully}. For the classical cases of $\alpha = 1$ and $\alpha = 2$, numerous numerical methods have been developed, such as finite difference methods \cite{dwivedi2023convergence,holden2015convergence,wang2021high,courtes2020error,thomee1998numerical}, Galerkin methods \cite{dutta2015convergence, dutta2016note,galtung2018convergent}, discontinuous Galerkin methods \cite{yan2002local,dwivedi4943062fully}, and spectral methods \cite{baker1983convergence,deng2009optimal,maday1988error,pelloni2001error}. 

\subsection{Challenges and Motivations}

While numerical schemes for the KdV and BO equations have been well studied, challenges remain in the convergence analysis of the Fourier Galerkin scheme to the unique solution of the corresponding equations. Specifically, for spectral methods, a significant gap exists in establishing rigorous convergence proofs for the approximations, particularly in non-local, fractional settings. For instance, the convergence of spectral Galerkin methods to the unique solution of the KdV and BO equations is still underdeveloped. Furthermore, the numerical approximation of the zero dispersion limit where $\varepsilon \to 0$ poses additional challenges, especially after the critical time $t_c$ when gradient catastrophe occurs, leading to the emergence of oscillatory behavior in the solution.
The foundational works of Lax and Levermore \cite{lax1983small, lax1983small2, lax1983small3} on the zero dispersion limit of the KdV equation, and numerical studies followed by Grava and Klein \cite{grava2012numerical,grava2007numerical},  established that numerical methods can accurately capture the oscillatory structures arising due to small dispersion coefficients. However, up to our knowledge, there is quite limited literature available on numerical schemes that effectively capture these oscillations, particularly for fractional KdV equations. This gap motivates our work to approximate the zero-dispersion limit of the fractional KdV equation \eqref{fkdv} accurately and {capture} the correct asymptotic behavior of solutions.

The spectral Galerkin method \cite{baker1983convergence,zayernouri2024spectral}, in particular, offers distinct advantages over traditional finite difference or finite element approaches when dealing with periodic and smooth initial data. These methods are higher order accurate and possess efficient resolution of fine-scale structures due to their global basis functions, making them well-suited for handling the oscillatory nature of dispersive waves. However, the lack of comprehensive convergence analysis, especially for the fractional equations, necessitates a detailed investigation into their stability and convergence behavior. Moreover, we would like to study the spectral Galerkin method in the small dispersion limit.

\subsection{Objectives and Contributions}

The primary objective of this paper is to develop and analyze a fully discrete Crank-Nicolson Fourier spectral Galerkin (FSG) scheme for approximating solutions to the fractional KdV equation \eqref{fkdv} with periodic initial data. Our main contributions are as follows:

\begin{itemize}
    \item \textbf{Conservation and stability:} We demonstrate that the semi-discrete FSG scheme preserves the first three integral invariants for all $\alpha \in [1,2]$, ensuring structure preservation. For the fully discrete scheme, we establish $L^2$-conservation, which guarantees $L^2$-stability of the numerical solution. In addition, we prove that the numerical approximations are uniformly bounded in the periodic $H_p^{1+\alpha}$ space, and that the temporal derivative of the approximations is bounded in the periodic $L^2$ space, contributing to the overall stability and reliability of the scheme.

    \item \textbf{Convergence to unique solution:} Using compactness arguments, we prove that the numerical approximation converges to the unique solution of the fractional KdV equation \eqref{fkdv} in $C([0,T]; H_p^{1+\alpha}(\R))$. While our convergence proof, as by product, can be viewed as a constructive proof of the existence and uniqueness of the solution to the fractional KdV equation \eqref{fkdv}.

    \item \textbf{Error analysis:} We derive error estimates for the proposed scheme, demonstrating spectral accuracy in space and second order accuracy in time for periodic initial data in $H_p^r$ with $r \geq 1+\alpha$, and exponential accuracy for analytic initial data.

    \item \textbf{Zero dispersion limit:} We numerically analyze the zero dispersion limit of the fractional KdV equation \eqref{fkdv} by analyzing the behavior of solutions as $\varepsilon \to 0$. Prior to the critical time $t_c$, the solution of the classical KdV equation converges to that of the Hopf equation for smooth initial data as $\varepsilon$ approaches zero. We extend this result to the fractional KdV equation \eqref{fkdv} by using the Kato’s theory \cite{kato1983cauchy}. Beyond the critical time $t_c$, within the oscillation zone, the solution is locally described by Whitham’s averaged equation, while outside this zone, it aligns with the Hopf equation. Our numerical results indicate that the fully discrete FSG scheme captures this asymptotic behavior, showing convergence to the asymptotic solution of the classical KdV equation for values of $\alpha$ close to 2.

    \item \textbf{Numerical validation:} We provide numerical examples to validate the theoretical analysis and to demonstrate the convergence and accuracy of the proposed scheme in approximating the fractional KdV equation \eqref{fkdv}.

\end{itemize}
  {In this paper, $C$ denotes a generic constant whose value can change in each step and it is independent of both the spatial discretization parameter $N$ and time discretization parameter $\Delta t$.}
\subsection{Outline}

The remainder of the paper is structured as follows: In Section \ref{sec2}, we present key results on Fourier analysis and the fractional Laplacian instrumental for our subsequent analysis. Section \ref{sec3} introduces the semi-discrete FSG scheme, establishes its conservation properties, and further, we provide the stability and convergence results for the fully discrete scheme, along with the error analysis. Section \ref{sec4} investigates the zero-dispersion limit of the fractional KdV equation. Numerical examples validating our theoretical findings are presented in Section \ref{sec5}. Finally, Section \ref{sec6} presents possible extensions of this work and provides concluding remarks.

\section{Notations and Preliminary Results}\label{sec2}

\subsection{Periodic Sobolev Spaces and Fourier Analysis}
We consider functions defined on $\mathbb{R}$ that are $2\pi$-periodic. For simplicity, we restrict these functions to the interval $I := [-\pi, \pi]$ and extend them periodically to the entire real line. The Hilbert space $L^2(I)$ consists of all square-integrable functions over $I$, with the norm and inner product defined by
\[
\norm{u} = \left(\int_I |u(x)|^2 \, dx\right)^{1/2} \quad \text{and} \quad (u,v) = \int_I u(x)\overline{v(x)} \, dx
\]
respectively, for all $u, v \in L^2(I)$, where $\overline{v}$ denotes the complex conjugate of $v$.

For $r \geq 0$, the Sobolev space $H^r(I) := W^{r,2}(I)$ is equipped with the norm $\norm{\cdot}_r$. We denote the periodic Sobolev spaces of exponent $r \geq 0$ by $H_p^r(I)$, consisting of all $2\pi$-periodic functions in $H^r(I)$. For a function $f \in H_p^r(I)$, the norm is defined as
\[
\norm{f}_r := \left(\sum_{k=-\infty}^{\infty}(1 + |k|^2)^r |\hat{f}(k)|^2\right)^{1/2},
\]
where $\hat{f}$ represents the Fourier coefficient of $f$, given by
\[
\hat{f}(k) := \frac{1}{2\pi}\int_I f(x)e^{-ikx} \, dx, \quad \forall k \in \mathbb{Z}.
\]
The Fourier expansion of $f \in H_p^r(I)$, $r \geq 0$, is given by
\begin{equation}\label{eqn_fourier}
f(x) = \sum_{k=-\infty}^{\infty} \hat{f}(k) e^{ikx},
\end{equation}
which converges almost everywhere. Clearly, $H_p^r(I)$ is a subspace of $H^r(I)$ with the norm $\|\cdot\|_r$, and $L_p^2(I) := H^0_p(I)$ with the norm $\|\cdot\|: = \|\cdot\|_0$.
 {The function space \( C_c^1\big([0, T]; H^{r}_p(I)\big) \) for \( r \geq 0 \) is defined as  
\begin{equation*}
{C_c^1\big([0, T]; H^{r}_p(I)\big) = \left\{ \varphi \in C^1([0,T]; H^r_p(I)) \;\middle|\;  \, 
\text{Supp}(\varphi(x,\cdot))\subset  [0, T] \text{ for fixed } x\in I \right\}.}
\end{equation*}

For any $N \in \mathbb{N}$, we define the approximation space consisting of real-valued trigonometric polynomials of degree $N$ by
\[
V_N := \text{span}\{e^{ikx} : -N \leq k \leq N\}.
\]
The spaces $V_N$, $N \in \mathbb{N}$, serve as natural approximations of $L_p^2(I)$ by periodic functions. Moreover, the basis functions $\{e^{ikx}\}_{k=-N}^{N}$ of $V_N$ are orthogonal. We define the projection operator $P_N: L^2(I) \to V_N$ by
\begin{equation}\label{eqn_projdefn}
P_N f(x) := \sum_{k=-N}^{N} \hat{f}(k) e^{ikx},
\end{equation}
and due to the orthogonality of the basis functions, $P_N$ is an orthogonal projection onto $V_N$. Equivalently, we have
\begin{equation}\label{eqn_ortho}
(P_N f - f, \phi) = 0, \quad \forall \, \phi \in V_N.
\end{equation}
Additionally, for any $f \in H^r(I)$, $r \geq 0$, it holds that $P_N f \to f$ in $L^2(I)$ as $N \to \infty$. In particular, we have the following result.

\begin{proposition}\label{Prop_Nr}
Let $r \geq 0$ be any real number and $r \geq s \geq 0$. Then, for $f \in H^r(I)$, there exists a positive constant $C$, independent of $N$, such that the following estimate holds
\begin{equation}\label{eqn_Nr}
\norm{f - P_N f}_s \leq C N^{-r+s} \norm{f}_r.
\end{equation}
\end{proposition}

The above result is well-known and can be found in \cite{hesthaven2017numerical}. However, for completeness, we provide a proof for the case $s = 0$.\\
\emph{Proof of Proposition \ref{Prop_Nr}}: Using Parseval's identity and the projection property \eqref{eqn_ortho}, we have
\begin{align*}
\norm{f - P_N f}^2 &= \norm{f}^2 - (f, P_N f) - (P_N f, f) + \norm{P_N f}^2 = 2\pi \sum_{k \in \mathbb{Z}} |\hat{f}(k)|^2 - 2\pi \sum_{k=-N}^{N} |\hat{f}(k)|^2 \\
&= 2\pi \sum_{|k| > N} |\hat{f}(k)|^2 = 2\pi \sum_{|k| > N} \frac{1}{(1 + |k|^2)^r} (1 + |k|^2)^r |\hat{f}(k)|^2 \\
&\leq 2\pi N^{-2r} \sum_{|k| > N} (1 + |k|^2)^r |\hat{f}(k)|^2 \leq C N^{-2r} \norm{f}_r^2.
\end{align*}
This completes the proof.\\

\begin{proposition}\label{Prop_exp}
Let $f$ be an analytic function. Then, there exist positive constants $C$ and $c$, independent of $N$, such that
\begin{equation}\label{eqn_eN}
\norm{f - P_N f} \leq C e^{-cN} \norm{f}.
\end{equation}
\end{proposition}
\begin{proof}
For all $r \geq 0$, it is known \cite{tadmor1986exponential} that
\[
\norm{f}_r \leq C r! \norm{f}.
\]
Utilizing Proposition \ref{Prop_Nr} and Stirling's formula \cite{tadmor1986exponential}, we have
\begin{align*}
\norm{f - P_N f} \leq C N^{-r} \norm{f}_r \leq C \frac{r!}{N^r} \norm{f} \leq C \frac{r^r e^{-r}}{N^r} \norm{f} \leq C e^{-cN} \norm{f},
\end{align*}
where we have assumed that $r$ is proportional to $N$. This completes the proof.
\end{proof}

\subsection{Fractional Laplacian and its properties}
We define the fractional Laplacian $\mathcal{D}^{\alpha}$ for $\alpha \geq 0$ using the Fourier series expansion \eqref{eqn_fourier} of a function $f \in H^r_p(I)$ as follows
\begin{equation}\label{eqn_fracLap}
    \mathcal{D}^{\alpha}f(x) = \sum\limits_{k\in\Z}|k|^\alpha \hat{f}(k) e^{ikx}.
\end{equation}
The following lemma states several important properties of the fractional Laplacian \eqref{eqn_fracLap}.

\begin{lemma}
    The fractional Laplacian \eqref{eqn_fracLap} satisfies the following properties
    \begin{enumerate}[label=\roman*)]
        \item For $f, g \in H^{\alpha}_p(I)$, $\alpha \geq 0$, we have
        \begin{equation}\label{eqn_fracsymm}
            \left(\mathcal{D}^{\alpha}f, g\right) = \left(f, \mathcal{D}^{\alpha}g\right),
        \end{equation}
        and
        \begin{equation}\label{eqn_fracortho}
            \left(\mathcal{D}^{\alpha}f_x, f\right) = 0.
        \end{equation}

        \item  Let $\alpha_1, \alpha_2 \geq 0$, then for all $f, g \in H_p^{\alpha_1+\alpha_2}(I)$, we have
        \begin{equation}\label{eqn_semigp}
            \left(\mathcal{D}^{\alpha_1+\alpha_2}f, g\right) = \left(\mathcal{D}^{\alpha_1}f, \mathcal{D}^{\alpha_2}g\right),
        \end{equation}
        and equivalently,
        \begin{equation}\label{eqn_semigp2}
            \mathcal{D}^{\alpha_1+\alpha_2}f = \mathcal{D}^{\alpha_1}\mathcal{D}^{\alpha_2}f = \mathcal{D}^{\alpha_2}\mathcal{D}^{\alpha_1}f.
        \end{equation}

        \item For an orthogonal projection $P_N$ defined by \eqref{eqn_ortho} and $f \in H^r_p(I)$, $r \geq \alpha \geq 0$, the fractional Laplacian \eqref{eqn_fracLap} with exponent $\alpha$ commutes with $P_N$, i.e.,
        \begin{equation}\label{eqn_projcom}
            \mathcal{D}^{\alpha} (P_N f(x)) = P_N \mathcal{D}^{\alpha}f(x).
        \end{equation}
    \end{enumerate}
\end{lemma}

\begin{proof}
    By using the Fourier expansion \eqref{eqn_fourier} for $f$ and $g$, and the orthogonality of $\{e^{ikx}\}$, we have
    \begin{align*}
        \left(\mathcal{D}^{\alpha}f, g\right) &= \left(\sum\limits_{k\in\Z}|k|^{\alpha} \hat{f}(k) e^{ikx}, \sum\limits_{k\in\Z} \hat{g}(k) e^{ikx}\right) = \sum\limits_{k\in\Z}|k|^{\alpha} \hat{f}(k) \hat{g}(k) \\&= \left(\sum\limits_{k\in\Z} \hat{f}(k) e^{ikx}, \sum\limits_{k\in\Z} |k|^{\alpha}\hat{g}(k) e^{ikx}\right) = \left(f, \mathcal{D}^{\alpha}g\right).
    \end{align*}
    {Now using \eqref{eqn_fracsymm} and periodicity of $f$, we have 
    $$\left(\mathcal{D}^{\alpha}f_x, f\right) = -\left(\mathcal{D}^{\alpha}f, f_x\right) =-\left(f,\mathcal{D}^{\alpha}f_x\right)=-\left(\mathcal{D}^{\alpha}f_x, f\right),$$
  which yields \eqref{eqn_fracortho}.} Similarly, we have
   \begin{align*}
        \left(\mathcal{D}^{\alpha_1+\alpha_2}f,g\right)&= \left( \sum\limits_{k\in\Z}|k|^{\alpha_1+\alpha_2} \hat f(k)e^{ikx},\sum\limits_{k\in\Z} \hat g(k)e^{ikx}\right)= \sum\limits_{k\in\Z}|k|^{\alpha_1+\alpha_2} \hat f(k) \hat g(k)\\& =  \sum\limits_{k\in\Z}|k|^{\alpha_1}|k|^{\alpha_2} \hat f(k) \hat g(k)= \left( \sum\limits_{k\in\Z}|k|^{\alpha_1}\hat f(k)e^{ikx},\sum\limits_{k\in\Z}|k|^{\alpha_2} \hat g(k)e^{ikx}\right)\\
        & = \left(\mathcal{D}^{\alpha_1}f, \mathcal{D}^{\alpha_2}g\right).
    \end{align*}
    which implies \eqref{eqn_semigp}. The identity \eqref{eqn_semigp2} follows directly from the definition \eqref{eqn_fracLap}.

    For \eqref{eqn_projcom}, using definitions \eqref{eqn_projdefn} and \eqref{eqn_fracLap}, we obtain
    \begin{align*}
        \mathcal{D}^{\alpha} (P_N f(x)) &= \sum\limits_{k=-N}^N |k|^{\alpha} \hat{f}(k) e^{ikx} = P_N \mathcal{D}^{\alpha}f(x).
    \end{align*}
    This completes the proof.
\end{proof}

\begin{lemma}\label{lemma_comest}
    Let $f, g \in H^{\alpha}_p(I)$, $\alpha \geq 0$. Then the following estimate holds
    \begin{equation}\label{eqn_lemma_comest}
        \mathcal{D}^\alpha(fg) \leq C(\alpha) \left(f \mathcal{D}^\alpha g + g \mathcal{D}^\alpha f\right),
    \end{equation}
    where $C(\alpha)$ is a constant depending on $\alpha$. Furthermore, there holds
    \begin{equation}\label{eqn_lemma_comest_morm}
        \norm{\mathcal{D}^\alpha(fg)} \leq C(\alpha) \left(\norm{f}_\infty \norm{\mathcal{D}^\alpha g} + \norm{g}_\infty \norm{\mathcal{D}^\alpha f}\right),
    \end{equation}
     {where $C(\alpha) = \max \left\{ 1, 2^{\alpha - 1} \right\}.$}
\end{lemma}

\begin{proof}
    To prove \eqref{eqn_lemma_comest}, we start by considering the Fourier series representation of the product $fg$. The Fourier coefficient of $fg$ is given by
    \begin{equation}\label{eqn_convl}
        \widehat{fg}(m) = \sum\limits_{k\in\Z} \hat{f}(k) \hat{g}(m-k).
    \end{equation}
    Applying the fractional derivative operator $\mathcal{D}^\alpha$ to the product $fg$, we obtain
    \begin{align*}
        \mathcal{D}^\alpha(fg)(x) &= \sum\limits_{m\in\Z} |m|^\alpha \widehat{fg}(m) e^{imx} = \sum\limits_{m\in\Z} |m|^\alpha \left( \sum\limits_{k\in\Z} \hat{f}(k) \hat{g}(m-k) \right) e^{imx}.
    \end{align*}
    We employ the inequality $|m|^\alpha \leq C(\alpha) \left(|k|^\alpha + |m-k|^\alpha\right)$ for $\alpha \geq 0$ where $m, k \in \Z$ and $C(\alpha) = \max \left\{ 1, 2^{\alpha - 1} \right\}$. Furthermore, incorporating the inequality \eqref{eqn_convl}, we can estimate the expression as follows
     {\begin{align*}
        \mathcal{D}^\alpha(fg)(x) &\leq C(\alpha) \Big[\sum\limits_{m\in\Z}    \sum\limits_{k\in\Z} {|k|^\alpha}\hat f(k)\hat g(m-k)e^{imx}\\  \nonumber      &\qquad+\sum\limits_{m\in\Z}   \sum\limits_{k\in\Z} {|m-k|^\alpha}\hat f(k)\hat g(m-k)e^{imx}\Big]\\ \nonumber
        &\leq C(\alpha) \Big[\sum\limits_{m\in\Z}\sum\limits_{k\in\Z}  \widehat{\mathcal{D}^\alpha f}(k) \hat g(m-k) e^{imx} +\sum\limits_{m\in\Z}\sum\limits_{k\in\Z}  \hat f(k) \widehat{\mathcal{D}^\alpha g}(m-k)e^{imx}\Big]\\
        &\leq C(\alpha) \Big[\sum\limits_{m\in\Z}  \widehat{g\mathcal{D}^\alpha f}(m) e^{imx} +\sum\limits_{m\in\Z}  \widehat{f\mathcal{D}^\alpha g}(m)e^{imx}\Big]\\
        &=C(\alpha) \big(g(x)\mathcal{D}^\alpha f(x)+ f(x)\mathcal{D}^\alpha g(x)\big).
    \end{align*}}
    This proves inequality \eqref{eqn_lemma_comest}.
    To prove \eqref{eqn_lemma_comest_morm}, we take the $L^2$-norm on both sides of \eqref{eqn_lemma_comest}
    \begin{align*}
        \norm{\mathcal{D}^\alpha(fg)} &\leq C(\alpha) \left( \norm{f}_\infty \norm{\mathcal{D}^\alpha g} + \norm{g}_\infty \norm{\mathcal{D}^\alpha f} \right).
    \end{align*}
Hence the result follows.
\end{proof}


\section{Fourier-Spectral-Galerkin Scheme: Stability and Convergence} \label{sec3}
The Fourier-Spectral-Galerkin (FSG) scheme for the fractional KdV equation \eqref{fkdv} is formulated as follows. Let $u_0\in H^{1+\alpha}_p(I)$, we seek an approximation $ U_\varepsilon \in V_N $ such that for all $\phi \in V_N$, the approximation $ U_\varepsilon $ satisfies
\begin{align}\label{eqn_fsgscheme}
    \begin{cases}
        \left((U_\varepsilon)_t + 6U_\varepsilon (U_\varepsilon)_x - \varepsilon^2 \mathcal{D}^{\alpha} (U_\varepsilon)_x, \phi\right) = 0, \qquad 0 \leq t \leq T, \ \varepsilon \in \mathbb{R}, \\
        U_\varepsilon(0) = P_N u_0,
    \end{cases}
\end{align}
where $ P_N $ denotes the projection operator and is defined by \eqref{eqn_projdefn}. Given the nonlinearity present in the scheme, it is crucial to establish the existence and uniqueness of the solution to \eqref{eqn_fsgscheme}. We state this in the following lemma.  {The approximation \(U_\varepsilon \in V_N\) represents the spectral Galerkin projection of the exact solution to \eqref{fkdv}, where the subscript emphasizes the dependence on the dispersion coefficient $\varepsilon$. Following the work of Kenig et al. \cite{kenig1991well}, the fractional KdV equation \eqref{fkdv} possesses the following conserved quantities:
\begin{align*}
    I_1(u):=\int_{I} u \, dx,\qquad
    I_2(u):= \int_{I} u^2 \, dx,\qquad
    I_3(u):= \int_{I} \left(\varepsilon^2(\mathcal{D}^{\alpha/2}u)^2 - 2u^3\right) dx,
\end{align*}
where $I_1$, $I_2$ and $I_3$ represents the {\em mass, momentum and energy}, respectively. In the following lemma, we establish discrete analogs of these conserved quantities for our numerical scheme \eqref{eqn_fsgscheme}.}

\begin{lemma}\label{lem_exis}
    There exists a unique solution $ U_\varepsilon $ to the equation \eqref{eqn_fsgscheme}. Moreover, the solution satisfies the first three conserved quantities of the fractional KdV equation \eqref{fkdv}, specifically
    \begin{align}
        \label{eqn_1stint} \frac{\partial}{\partial t}\left[ \int_{-\pi}^{\pi} U_\varepsilon(x,t) \, dx \right] &= 0, \\
        \label{eqn_2ndint} \frac{\partial}{\partial t}\left[ \int_{-\pi}^{\pi} (U_\varepsilon(x,t))^2 \, dx \right] &= 0, \\
        \label{eqn_3rdint} \frac{\partial}{\partial t}\left[ \int_{-\pi}^\pi \left(\varepsilon^2 \left(\mathcal{D}^{\alpha/2} U_\varepsilon(x,t)\right)^2 - 2U^3_\varepsilon(x,t)\right) \, dx \right] &= 0.
    \end{align}   
\end{lemma}
\begin{proof}
    We choose test functions $\phi = e^{ikx}$, for \(k = -N, -N+1, \ldots, N-1, N\) in the scheme \eqref{eqn_fsgscheme} and use the fact that since \(u_0\) is real-valued, \(U_\varepsilon(0)\) is also real-valued. This implies \(\hat{U}_\varepsilon(k,t) = \overline{\hat{U}_\varepsilon(k,t)}\), leading to the following system of equations for the Fourier coefficients \(\hat{U}_\varepsilon(k,t)\) of \(U_\varepsilon\)
    \begin{align}\label{eqn_ft}
     \begin{cases}
        (\hat{U}_\varepsilon(k,t))_t  = -3ik \hat{U}_\varepsilon \ast \hat{U}_\varepsilon(k,t) + ik|k|^\alpha \varepsilon^2 \hat{U}_\varepsilon(k,t), \\
        \hat{U}_\varepsilon(k,0)  = \hat{u}_0(k), \end{cases}
    \end{align}
    for \(k = -N, -N+1, \ldots, N-1, N\). 
 {The right-hand side of \eqref{eqn_ft} is locally Lipschitz continuous in $\hat{U}_{\varepsilon}$ with respect to the \(L^2\) norm. By the Picard-Lindelöf Theorem, there exists a unique solution to \eqref{eqn_ft}. This implies that there exists a maximal time \(t_0, 0<t_0\leq T\), such that for all \(t<t_0\), \eqref{eqn_ft} admits a local unique solution, which implies \eqref{eqn_fsgscheme} has a unique solution \({U}_\varepsilon(t)\) for all $t<t_0$.} Now, choosing the test function \(\phi = U_\varepsilon\) in \eqref{eqn_fsgscheme} and using the property \eqref{eqn_fracortho} yields
    \begin{align*}
        ((U_\varepsilon)_t, U_\varepsilon) &= - 6(U_\varepsilon (U_\varepsilon)_x , U_\varepsilon) +\varepsilon^2 ( \mathcal{D}^{\alpha} (U_\varepsilon)_x, U_\varepsilon) = 0,
    \end{align*}
    which further implies
    \begin{equation}\label{eqn_L2c}
        \|U_\varepsilon(T)\| = \|U_\varepsilon(0)\| = \|u_0\|,
    \end{equation}
      which implies \eqref{eqn_2ndint} and ensures that the solution does not blow up for all \(t \leq T\), and \(T\) can be chosen sufficiently large. This establishes the global existence of the solution to \eqref{eqn_fsgscheme}. 
    To show \eqref{eqn_1stint}, we choose the test function \(\phi = 1\) in \eqref{eqn_fsgscheme} and use the periodicity of \(U_\varepsilon\), which gives \eqref{eqn_1stint}.
    In order to show \eqref{eqn_3rdint}, we choose the test function \(\phi = P_N\big(3U_\varepsilon^2 - \varepsilon^2 \mathcal{D}^{\alpha}U_\varepsilon\big)\) in \eqref{eqn_fsgscheme} to get
    \begin{align}\label{eqn_con1}
        \int_{-\pi}^\pi (U_\varepsilon)_t P_N\big(3U_\varepsilon^2 - \varepsilon^2 \mathcal{D}^{\alpha}U_\varepsilon\big)\,dx = - \int_{-\pi}^\pi \left(3U_\varepsilon^2 - \varepsilon^2 \mathcal{D}^{\alpha}U_\varepsilon\right)_x P_N\big(3U_\varepsilon^2 - \varepsilon^2 \mathcal{D}^{\alpha}U_\varepsilon\big)\,dx. 
    \end{align}
    We estimate the left-hand side of the above equation utilizing the identity \eqref{eqn_ortho} and using the fact that \((U_\varepsilon)_t\) is in \(V_N\), which implies
    \begin{align}\label{eqn_con2}
        \nonumber\int_{-\pi}^\pi (U_\varepsilon)_t P_N\big(3U_\varepsilon^2 - \varepsilon^2 \mathcal{D}^{\alpha}U_\varepsilon\big)\,dx &= \int_{-\pi}^\pi (U_\varepsilon)_t \big(3U_\varepsilon^2 - \varepsilon^2 \mathcal{D}^{\alpha}U_\varepsilon\big)\,dx\\ 
        \nonumber &= \int_{-\pi}^\pi(U_\varepsilon^3)_t - \varepsilon^2 (U_\varepsilon)_t \mathcal{D}^{\alpha}U_\varepsilon\,dx\\
        &= \frac{\partial}{\partial t}\int_{-\pi}^\pi \left(U^3_\varepsilon(x,t) -\frac{1}{2} \varepsilon^2 \big(\mathcal{D}^{\alpha/2}U_\varepsilon(x,t)\big)^2\right)\,dx,
    \end{align}
    where we have used the semigroup property \eqref{eqn_semigp} of the fractional Laplacian.
    On the other hand, we estimate the right-hand side of \eqref{eqn_con1} by utilizing the orthogonal property \eqref{eqn_ortho} of \(P_N\) and using the fact that \(P_N\) commutes with derivatives
    \begin{align}\label{eqn_con3}
       \int_{-\pi}^\pi\left( 3U_\varepsilon^2 - \varepsilon^2 \mathcal{D}^{\alpha}(U_\varepsilon)\right)_x P_N\big(3U_\varepsilon^2 - \varepsilon^2 \mathcal{D}^{\alpha}U_\varepsilon\big)\,dx = \frac{1}{2} \int_{-\pi}^\pi \left(P_N\big(3U_\varepsilon^2 - \varepsilon^2 \mathcal{D}^{\alpha}U_\varepsilon\big)\right)^2_x \,dx
        = 0.
    \end{align}
    Thus, from identities \eqref{eqn_con1}, \eqref{eqn_con2}, and \eqref{eqn_con3}, we obtain \eqref{eqn_3rdint}. This completes the proof.
\end{proof}

We employ the Crank-Nicolson method for time discretization in \eqref{eqn_fsgscheme} to define the fully discrete FSG scheme. Let \(\Delta t\) be the time step size, and define \(U_{\varepsilon}^m = U_{\varepsilon}(t_m)\), where \(t_m = m\Delta t\) for \(m = 0, 1, \ldots, M\), with \(t_M = T\) for a given \(0 < T < \infty\). The fully discrete FSG scheme is formulated as follows: given \(U_{\varepsilon}^m \in V_N\), find \(U^{m+1}_{\varepsilon} \in V_N\) such that, for all test functions \(\phi \in V_N\), the scheme satisfies
\begin{align}\label{eqn_fdfsgscheme}
    \begin{cases}
        \left(U^{m+1}_\varepsilon , \phi\right) = \left(U^{m}_\varepsilon, \phi\right) 
        - 6\Delta t\left(U^{m+\frac{1}{2}}_\varepsilon (U^{m+\frac{1}{2}}_\varepsilon)_x , \phi\right) 
        + \varepsilon^2 \Delta t \left(\mathcal{D}^{\alpha} (U^{m+\frac{1}{2}}_\varepsilon)_x, \phi\right), \\
        U^0_\varepsilon = P_N u_0,
    \end{cases}
\end{align}
for \(m = 0, 1, \ldots, M-1\) and \(U^{m+\frac{1}{2}}_\varepsilon = (U^m_\varepsilon + U^{m+1}_\varepsilon)/2\). Clearly, we have \(\|U^0_\varepsilon\| = \|P_N u_0\| \leq \|u_0\|\).
Note that the fully discrete scheme \eqref{eqn_fdfsgscheme} conserves the \(L^2\) norm, leading to the following stability result.
\begin{lemma}\label{lemma_stab}
    The fully discrete scheme \eqref{eqn_fdfsgscheme} is $L^2-stable$. Moreover, for the approximation $U^m_{\varepsilon}$, we have the following $L^2$-bound 
    \begin{equation}\label{eqn_L2stabfd}
        \norm{U_{\varepsilon}^m} \leq C, \qquad \text{ for all } n \text{ and } \varepsilon,
    \end{equation}
    where $C$ is a constant independent of $m$ and $\varepsilon$.
\end{lemma}
\begin{proof}
    By choosing the test function $\phi = U^{m+\frac{1}{2}}_{\varepsilon}$ in \eqref{eqn_fdfsgscheme} and utilizing identity \eqref{eqn_fracortho} yields
    \begin{equation*}
        \norm{U_\varepsilon^{m+1}}^2 = \norm{U_{\varepsilon}^m}^2,
    \end{equation*}
    which implies the estimate \eqref{eqn_L2stabfd} and scheme \eqref{eqn_fdfsgscheme} is $L^2$-stable.
\end{proof}

\subsection{Solvability of the fully discrete FSG scheme}

Given that the non-linear part of the scheme involves an implicit term, it is necessary to address the solvability of the system at each time step. To solve the non-linear system, we introduce an iterative sequence \(\{v^\ell\}_{\ell \geq 0}\) such that for each \(\ell \geq 0\), \(v^\ell \in V_N\) satisfy the following system
\begin{equation}\label{eqn_itera}
    \begin{cases}
        (v^{\ell+1}, \phi) = (U^m_{\varepsilon}, \phi) 
        - \Delta t \left( \mathcal{B}\left(\frac{U^m_\varepsilon + v^\ell}{2}\right), \phi \right) 
        + \frac{1}{2}\varepsilon^2 \Delta t\left( \mathcal{D}^\alpha (U^m_\varepsilon + v^{\ell+1})_x, \phi \right), \\
        v^0 = U^m_{\varepsilon},
    \end{cases}
\end{equation}
for all test functions \(\phi \in V_N\), where \(\mathcal{B}(v) = 6vv_x\).
 {The above iterative scheme can be considered as an FSG scheme for a linear problem involving the bilinear form in \(v^{\ell+1}\) and \(U^m_\varepsilon\). The scheme \eqref{eqn_itera} can be rewritten as  
\begin{equation*}\label{solv2}
    \left( \left(1-\frac{\varepsilon^2\Delta t}{2}\mathcal{D}^\alpha\partial_x \right)v^{\ell+1},\phi \right) =(U^m_{\varepsilon}, \phi) 
        - \Delta t \left( \mathcal{B}\left(\frac{U^m_\varepsilon + v^\ell}{2}\right), \phi \right) 
        + \frac{1}{2}\varepsilon^2 \Delta t\left( \mathcal{D}^\alpha(U^m_\varepsilon)_x, \phi \right).
\end{equation*}
Setting \(d_{ij} = \left( \mathcal{D}^\alpha\partial_x \phi_j,\phi_i \right) = - \left(\mathcal{D}^\alpha\partial_x \phi_i,\phi_j \right) = - d_{ji}\), where \(\phi_j\in V_N\) for all \(-N\leq i,j\leq N\), forming an orthogonal basis for \(V_N\). Choosing \(\phi = \phi_j\) for \(-N\leq j\leq N\) and writing $v^{\ell+1}$ as a linear combination of $\phi_j$, the left-hand side coefficient matrix becomes \(I-\frac{\varepsilon^2\Delta t}{2}D\), where \(D=(d_{ij})\) is skew-symmetric. Since \(I-\frac{\varepsilon^2\Delta t}{2}D\) is non-singular, the existence and uniqueness of \(v^{\ell+1}\) in \eqref{eqn_itera} is ensured.}
The solvability of the iterative scheme \eqref{eqn_itera} is established in the next lemma. Specifically, we prove that for a given \(U_{\varepsilon}^m\), the iteration $v^\ell$ converges to the next time step solution \(U_\varepsilon^{m+1}\) as $\ell\to\infty$.   Additionally, we prove that the fully discrete scheme \eqref{eqn_fdfsgscheme} remains stable in higher-order Sobolev spaces. We begin by demonstrating the stability of the approximation at each time step.

\begin{lemma}\label{solv_lemma}
 {Assume \( U^m_\varepsilon \) is an approximate} solution obtained by the fully discrete scheme \eqref{eqn_fdfsgscheme} at the \(m\)-th time step. Suppose \( (m+1)\Delta t \leq T \), for some $T>0$, and let \( v^\ell \) solve the iterative scheme \eqref{eqn_itera}. Furthermore, assume that \(\Delta t\) is sufficiently small such that at the \(m\)-th step it satisfies
\begin{equation}\label{eqn_timestep}
    N \Delta t \leq  \frac{\zeta}{C\eta \norm{U^m_\varepsilon}_{1+\alpha}},
\end{equation}
where \(\zeta \in (0,1)\), \(\eta = \frac{8-\zeta}{1-\zeta} > 8\), $\alpha\in[1,2]$, and  {$C$ is some constant independent of $\Delta t$ and $N$.} Then the sequence \(\{v^\ell\}_{\ell\geq0}\) converges, and
\begin{equation}
    \lim\limits_{\ell\to\infty} v^\ell = U^{m+1}_\varepsilon.
\end{equation}
Moreover, the approximation \( U^{m+1}_\varepsilon \)  satisfies
\begin{equation}
    \norm{U^{m+1}_\varepsilon}_{1+\alpha} \leq \eta \norm{U^m_\varepsilon}_{1+\alpha}.
\end{equation}
\end{lemma}

\begin{proof}
We start with the iteration \eqref{eqn_itera} for \(\ell=0\). For \(v^1\), it satisfies
\begin{equation}\label{eqn_aux1}
    (v^{1}, \phi) = (U^m_{\varepsilon}, \phi) - \Delta t \left( \mathcal{B}(U^m_\varepsilon), \phi \right) + \frac{\Delta t}{2}\varepsilon^2 \left( \mathcal{D}^\alpha (U^m_\varepsilon + v^{1})_x, \phi \right).
\end{equation}
Choosing the test function \(\phi  = \mathcal{D}^{2\alpha}(U^m_{\varepsilon}+v^1)_{xx} { \in V_N}\) in \eqref{eqn_aux1}, and using identity \eqref{eqn_fracortho} and Young's inequality, we have
\begin{align*}
    \nonumber\norm{\mathcal{D}^\alpha v^1_{x}}^2 &= \norm{\mathcal{D}^\alpha(U^m_{\varepsilon})_{x}}^2 +\Delta t \left(\mathcal{D}^\alpha\mathcal{B}(U^m_\varepsilon)_{x}, \mathcal{D}^\alpha(U^m_{\varepsilon}+v^1)_{x}\right)\\
       &\leq \norm{\mathcal{D}^\alpha(U^m_{\varepsilon})_{x}}^2 +  {2}\Delta t^2 \norm{\mathcal{D}^\alpha\mathcal{B}(U^m_\varepsilon)_{x}}^2 + \frac{1}{4} \norm{\mathcal{D}^\alpha(U^m_{\varepsilon})_{x}}^2 + \frac{1}{4} \norm{\mathcal{D}^\alpha v^1_{x}}^2,
\end{align*}
implies
\begin{align}\label{eqn_aux2}
   \norm{\mathcal{D}^\alpha v^1_{x}}^2 \leq 2\norm{\mathcal{D}^\alpha(U^m_{\varepsilon})_{x}}^2 + \frac{8}{3}\Delta t^2 \norm{\mathcal{D}^\alpha\mathcal{B}(U^m_\varepsilon)_{x}}^2.
\end{align}
 {Using the Lemma \ref{lemma_comest}, Lemma \ref{A.1} and $C(\alpha)=\max\{1,2^{\alpha-1}\}$, we obtain
\begin{align}\label{eqn_auxsob}
    \nonumber\norm{\mathcal{D}^\alpha\mathcal{B}(U^m_\varepsilon)_{x}}^2 &\leq 6\norm{\mathcal{D}^\alpha ((U^m_\varepsilon)_{x})^2}^2 +6\norm{ \mathcal{D}^\alpha (U^m_\varepsilon (U^m_\varepsilon)_{xx})}^2\\
&\nonumber\leq 6\Big(4C(\alpha)^2\|U^m_\varepsilon\|_{\infty}^2\|\mathcal{D}^\alpha (U^m_\varepsilon)_x\|^2+2C(\alpha)^2\|(U^m_\varepsilon)_{xx}\|_{\infty}^2\|\mathcal{D}^\alpha U^m_\varepsilon\|^2 \\ \nonumber \qquad &+2C(\alpha)^2\|U^m_\varepsilon\|_{\infty}^2\|\mathcal{D}^\alpha (U^m_\varepsilon)_{xx}\|^2\Big)\\
    &\leq CN^2 \norm{U^m_{\varepsilon}}_{1+\alpha}^4,
\end{align}}
for some constant $C$ independent of $\Delta t$ and $N$. Applying the Sobolev inequalities from Lemma \ref{A.1} and estimate \eqref{eqn_auxsob} in \eqref{eqn_aux2}, we get
\begin{align}\label{eqn_aux3}
    \norm{\mathcal{D}^\alpha v^1_{x}} \leq (2 +  {C}N^2\Delta t^2 \norm{U^m_{\varepsilon}}_{1+\alpha}^2)^{1/2} \norm{U^m_{\varepsilon}}_{1+\alpha}.
\end{align}
The time step condition \eqref{eqn_timestep} implies 
\begin{equation*}
    (2 +  {C}N^2\Delta t^2 \norm{U^m_{\varepsilon}}_{1+\alpha}^2)^{1/2} \leq 2.
\end{equation*}
As a consequence, \eqref{eqn_aux3} reduces to
\begin{align}\label{eqn_aux4}
    \norm{\mathcal{D}^\alpha v^1_{x}} \leq 2 \norm{U^m_{\varepsilon}}_{1+\alpha}.
\end{align}
Similarly, we estimate the lower order derivatives \(\norm{v^1_{xx}}, \norm{v^1_{x}}\), and \(\norm{v^1}\), which leads to
\begin{align}\label{eqn_aux5}
    \norm{v^1}_{1+\alpha} \leq \eta \norm{U^m_{\varepsilon}}_{1+\alpha}.
\end{align}
Afterwards, the iteration \eqref{eqn_itera} can be rewritten as
\begin{equation*}
    \left( \Big( 1 - \frac{\Delta t}{2}\varepsilon^2 \mathcal{D}^\alpha\partial_x\Big) v^{\ell+1}, \phi\right) = (U^m_{\varepsilon}, \phi) - \Delta t \left(\mathcal{B}\Big(\frac{U^m_\varepsilon + v^\ell}{2}\Big), \phi\right) + \frac{\Delta t}{2}\varepsilon^2 \left( \mathcal{D}^\alpha (U^m_\varepsilon)_x, \phi\right).
\end{equation*}
Subtracting the consecutive iterations, we have
\begin{equation}\label{eqn_aux6}
    \left( \Big( 1 - \frac{\Delta t}{2}\varepsilon^2 \mathcal{D}^\alpha\partial_x\Big) \delta (v^{\ell}), \phi\right) = - \Delta t \left(\mathcal{B}\Big(\frac{U^m_\varepsilon + v^\ell}{2}\Big) - \mathcal{B}\Big(\frac{U^m_\varepsilon + v^{\ell-1}}{2}\Big), \phi\right), 
\end{equation}
where \(\delta (v^\ell) = v^{\ell+1} - v^{\ell}\). Choosing \(\phi = \mathcal{D}^{2\alpha}\delta(v^\ell)_{xx}\in V_N\) in \eqref{eqn_aux6}, we have
 \begin{align*}
    \|\mathcal{D}^\alpha\delta (v^{\ell})_{x}\|^2 &= -\Delta t \left(\mathcal{D}^\alpha\mathcal{B}\Big(\frac{U^m_\varepsilon + v^\ell}{2}\Big)_{x}-\mathcal{D}^\alpha\mathcal{B}\Big(\frac{U^m_\varepsilon + v^{\ell-1}}{2}\Big)_{x}, \mathcal{D}^\alpha\delta (v^{\ell})_{x}\right)\\
     &\leq \Delta t \norm{\mathcal{D}^\alpha\mathcal{B}\Big(\frac{U^m_\varepsilon + v^\ell}{2}\Big)_{x}-\mathcal{D}^\alpha\mathcal{B}\Big(\frac{U^m_\varepsilon + v^{\ell-1}}{2}\Big)_{x}}\norm{\mathcal{D}^\alpha\delta (v^{\ell})_{x}}\\
     &\leq CN\Delta t \max\{\|U^m_\varepsilon\|_2, \|v^\ell\|_{1+\alpha},\|v^{\ell-1}\|_{1+\alpha}\}\|\delta(v^{\ell-1})\|_{1+\alpha}\|\mathcal{D}^\alpha\delta(v^{\ell-1})_{x}\|,
 \end{align*}
which simplifies to
\begin{equation*}
    \|\mathcal{D}^\alpha\delta (v^{\ell})_{x}\| \leq CN\Delta t \max\{\|U^m_\varepsilon\|_{1+\alpha}, \|v^\ell\|_{1+\alpha},\|v^{\ell-1}\|_{1+\alpha}\}\|\delta(v^{\ell-1})\|_{1+\alpha}.
\end{equation*}
Similar estimates can be obtained for the lower order derivatives as well. Consequently, we end up with
\begin{equation}\label{eqn_aux7}
    \|\delta (v^{\ell})\|_{1+\alpha} \leq CN\Delta t \max\{\|U^m_\varepsilon\|_{1+\alpha}, \|v^\ell\|_{1+\alpha}, \|v^{\ell-1}\|_{1+\alpha}\}\|\delta(v^{\ell-1})\|_{1+\alpha}.
\end{equation}
In particular, for $\ell = 1$, estimate \eqref{eqn_aux7} yields
\begin{align}\label{eqn_aux8}
    \|\delta (v^{1})\|_{1+\alpha} &\leq CN\Delta t \max\{\|U^m_\varepsilon\|_{1+\alpha}, \|v^1\|_{1+\alpha}\}\|\delta(v^{0})\|_{1+\alpha}\\
    & \leq CN\Delta t \eta \|U^m_\varepsilon\|_{1+\alpha} \|\delta(v^{0})\|_{1+\alpha} \nonumber \leq \zeta \|\delta(v^{0})\|_{1+\alpha},
\end{align}
where we have used the estimate \eqref{eqn_aux5} and time step condition \eqref{eqn_timestep}. With the help of estimates \eqref{eqn_aux5} and \eqref{eqn_aux8}, we proceed by induction. Let us assume that the following estimates hold
\begin{align}
    \label{eqn_aux9} \|v^\ell\|_{1+\alpha} &\leq \eta \norm{U^m_\varepsilon}_{1+\alpha},\\
    \label{eqn_aux10} \|\delta(v^\ell)\|_{1+\alpha} &\leq \zeta \|\delta(v^{\ell-1})\|_{1+\alpha},
\end{align}
for $\ell = 1, 2, \cdots, d$. We claim that the above estimates hold for $\ell = d+1$. Using the triangle inequality, assumption \eqref{eqn_aux10}, and time step condition \eqref{eqn_timestep}, we obtain
\begin{align*}
    \|v^{d+1}\|_{1+\alpha} &\leq \sum\limits_{\ell=0}^d \|\delta(v^\ell)\|_{1+\alpha} + \norm{U^m_\varepsilon}_{1+\alpha} \leq \|v^1 - U^m_\varepsilon\|_{1+\alpha} \sum\limits_{\ell=0}^d \zeta^\ell + \norm{U^m_\varepsilon}_{1+\alpha} \\ 
    &\leq \frac{1}{1-\zeta}(\|v^1\|_{1+\alpha} + \|U^m_\varepsilon\|_{1+\alpha}) + \norm{U^m_\varepsilon}_{1+\alpha} \leq \frac{8-\zeta}{1-\zeta} \norm{U^m_\varepsilon}_{1+\alpha} = \eta \norm{U^m_\varepsilon}_{1+\alpha}.
\end{align*}
Again using the time step condition \eqref{eqn_timestep} and estimate \eqref{eqn_aux7}, we have
\begin{equation}\label{eqn_aux12}
    \|\delta(v^{d+1})\|_{1+\alpha} \leq CN\Delta t \eta \|U^m_\varepsilon\|_{1+\alpha} \|\delta(v^{d})\|_{1+\alpha} \leq \zeta \|\delta(v^{d})\|_{1+\alpha}.
\end{equation}
Hence the estimate \eqref{eqn_aux10} is obtained for all $d \geq 0$. Moreover, estimate \eqref{eqn_aux10} implies that $v^\ell$ is a Cauchy sequence, and hence it converges, completing the proof of the claim.
\end{proof}

 \begin{remark}
    We have established that the devised scheme \eqref{eqn_fdfsgscheme} is solvable at time $(m+1)\Delta t$, given the approximate solution $U^m_\varepsilon$ at time $m\Delta t$, under the time step restriction \eqref{eqn_timestep}. Note that this condition depends on the $H^{1+\alpha}$-norm of $U^m_\varepsilon$, rather than directly on the initial data. However, to ensure that the CFL condition depends only on the initial data $u_0$, it is necessary to derive \emph{a priori} bounds for the numerical solution $U_{\varepsilon}^m$. This will be established in the next lemma.
\end{remark}

\begin{lemma}\label{lemma_pbound}
    Let $U_{\varepsilon}^m$ be an approximate solution to equation \eqref{fkdv} obtained via the fully discrete scheme \eqref{eqn_fdfsgscheme}. Assume that the initial data $u_0 \in H_p^{1+\alpha}(I)$, $\alpha \in [1,2]$. Further, suppose that the time step satisfies the CFL condition
    \begin{equation}\label{eqn_cfl}
         N \Delta t \leq  \frac{\zeta}{C\eta Z},
    \end{equation}
    where $Z$ is a constant depending on $\|u_0\|_{1+\alpha}$, and $\zeta$, $\eta$ and $C$ are as defined in Lemma \ref{solv_lemma}.  {Then there exists a finite time $\bar T$ such that} the following estimates hold:
    \begin{align}
        \label{eqn_Ubd3} \|U_{\varepsilon}^m\|_{1+\alpha} &\leq C, \\
        \label{eqn_Ubd2} \|D_+^t U_{\varepsilon}^m\| &\leq C,
    \end{align}
    for all $m \Delta t \leq \bar{T}$, where  {$D_+^t U_{\varepsilon}^m= \frac{1}{\Delta t}(U_\varepsilon^{m+1}- U_{\varepsilon}^m)$,} and $C$ depends only on $\alpha$, $\varepsilon$, and $\bar T$.
\end{lemma}

\begin{proof}
    We choose $\phi = \mathcal{D}^{2\alpha} (U^{m+\frac{1}{2}}_\varepsilon)_{xx}\in V_N$ as the test function in \eqref{eqn_fdfsgscheme}, yielding
    \begin{equation}\label{eqn_bddp1}
        \|\mathcal{D}^\alpha (U^{m+1}_\varepsilon)_x\|^2 - \|\mathcal{D}^\alpha (U^m_\varepsilon)_x\|^2 = 2 \Delta t \left( \mathcal{D}^\alpha \mathcal{B}(U^{m+\frac{1}{2}}_\varepsilon)_x, \mathcal{D}^\alpha (U^{m+\frac{1}{2}}_\varepsilon)_x \right),
    \end{equation}
    where we used the identity \eqref{eqn_fracortho} to eliminate the fractional term and periodicity to shift derivatives.
    To estimate the right-hand side of \eqref{eqn_bddp1}, we employ the  {Sobolev inequality from Lemma \ref{A.1}} and Lemma \ref{lemma_comest}. Applying multiple instances of the Cauchy-Schwarz inequality and Lemma \ref{lemma_comest}, we obtain
    \begin{small}
        \begin{align*}
        \Big|\Big( \mathcal{D}^\alpha \mathcal{B} &(U^{m+\frac{1}{2}}_\varepsilon)_x, \mathcal{D}^\alpha(U^{m+\frac{1}{2}}_\varepsilon)_{x} \Big)\Big|  \leq  6C(\alpha)\left|\left(U^{m+\frac{1}{2}}_\varepsilon \mathcal{D}^\alpha (U^{m+\frac{1}{2}}_\varepsilon)_{xx} + (U^{m+\frac{1}{2}}_\varepsilon)_{xx} \mathcal{D}^\alpha U^{m+\frac{1}{2}}_\varepsilon, \mathcal{D}^\alpha (U^{m+\frac{1}{2}}_\varepsilon)_x\right)\right| \\& \qquad \qquad+  12 C(\alpha) \|(U^{m+\frac{1}{2}}_\varepsilon)_x\|_\infty \|\mathcal{D}^\alpha (U^{m+\frac{1}{2}}_\varepsilon)_x\|^2 \\
        & \leq   3 C(\alpha) \left|\left((U^{m+\frac{1}{2}}_\varepsilon)_x, (\mathcal{D}^\alpha (U^{m+\frac{1}{2}}_\varepsilon)_x)^2\right)\right|  +6 C(\alpha) \|(U^{m+\frac{1}{2}}_\varepsilon)_{xx}\| \|\mathcal{D}^\alpha U^{m+\frac{1}{2}}_\varepsilon\|_\infty \|\mathcal{D}^\alpha (U^{m+\frac{1}{2}}_\varepsilon)_x\| \\ &\qquad\qquad + 12 C(\alpha) \|(U^{m+\frac{1}{2}}_\varepsilon)_x\|_\infty \|\mathcal{D}^\alpha (U^{m+\frac{1}{2}}_\varepsilon)_x\|^2 \\
        & \leq  {C} \|\mathcal{D}^\alpha (U^{m+\frac{1}{2}}_\varepsilon)_x\| \|U^{m+\frac{1}{2}}_\varepsilon\|_{1+\alpha}^2,
    \end{align*}
    \end{small}
     {for some constant $C$ independent of $\Delta t$ and $N$.}
    Substituting this into \eqref{eqn_bddp1} and using the triangle inequality, we get
    \begin{equation}\label{eqn_bddp2}
        \|\mathcal{D}^\alpha (U^{m+1}_\varepsilon)_x\| - \|\mathcal{D}^\alpha (U^m_\varepsilon)_x\| \leq 2 \Delta t \frac{C\|\mathcal{D}^\alpha (U^{m+\frac{1}{2}}_\varepsilon)_x\| \|U^{m+\frac{1}{2}}_\varepsilon\|_{1+\alpha}^2}{\|\mathcal{D}^\alpha (U^{m+1}_\varepsilon)_x\| + \|\mathcal{D}^\alpha (U^m_\varepsilon)_x\|} \leq C \Delta t \|U^{m+\frac{1}{2}}_\varepsilon\|_{1+\alpha}^2.
    \end{equation}
    A similar approach yields estimates for the terms with lower order derivatives, and combining with \eqref{eqn_bddp2}, we obtain
    \begin{equation}\label{eqn_bddp3}
        \|U^{m+1}_\varepsilon\|_{1+\alpha} \leq \|U^m_\varepsilon\|_{1+\alpha} +  C \Delta t \|U^{m+\frac{1}{2}}_\varepsilon\|_{1+\alpha}^2.
    \end{equation}
    Consider the differential equation
    \begin{equation}\label{eqn_odeaux}
        z'(t) =  \frac{C}{4} (1+\eta)^2 z(t)^2, \quad t > 0; \qquad z(0) = \|u_0\|_{1+\alpha},
    \end{equation}
    where $C$ is a constant independent of $t$. The solution $z(t)$ is increasing and convex for $t < T_b = \frac{4}{C (1+\eta)^2 z(0)}$. In particular, if $t_m \leq \bar T := \frac{T_b}{2}$, then $z(t_m) \leq z(\bar T) := Z$, and we aim to show
    \begin{equation}\label{eqn_help1}
        \|U^m_\varepsilon\|_{1+\alpha} \leq z(t_m) \leq Z, \quad \text{for } t_m \leq \bar T.
    \end{equation}
    The claim holds for $m = 0$. Assume it holds for $m = 1, 2, \cdots, d$. The CFL condition \eqref{eqn_cfl} ensures the time step condition \eqref{eqn_timestep}, and Lemma \ref{solv_lemma} provides
    \begin{equation}\label{eqn_help2}
        \|U^{m+\frac{1}{2}}_\varepsilon\|_{1+\alpha} \leq \frac{1+\eta}{2} \|U^m_\varepsilon\|_{1+\alpha}.
    \end{equation}
    Substituting this into \eqref{eqn_bddp3} and using the fact that $z$ is increasing and convex, we get
    \begin{align*}
        \|U^{d+1}_\varepsilon\|_{1+\alpha} &\leq \|U^d_\varepsilon\|_{1+\alpha} + \frac{C}{4} \Delta t (1+\eta)^2 \|U^d_\varepsilon\|_{1+\alpha}^2 \leq z(t_d) +  \frac{C}{4} \Delta t (1+\eta)^2 z(t_d)^2 \\
        &\leq z(t_d) + \int_{t_d}^{t_{d+1}}  \frac{C}{4} (1+\eta)^2 z(s)^2 \, ds
        \leq z(t_d) + \int_{t_d}^{t_{d+1}} z'(s) \, ds = z(t_{d+1}).
    \end{align*}
    This confirms \eqref{eqn_help1} and establishes \eqref{eqn_Ubd3}.
    To prove \eqref{eqn_Ubd2}, consider the scheme \eqref{eqn_fdfsgscheme} which implies
    \begin{equation}\label{eqn_help4}
        \left(D_+^t U^m_\varepsilon, \phi\right) = - \left(\mathcal{B}(U^{m+\frac{1}{2}}_\varepsilon), \phi\right) + \varepsilon^2 \left(\mathcal{D}^\alpha (U^{m+\frac{1}{2}}_\varepsilon)_x, \phi\right).
    \end{equation}
    Choosing $\phi = D_+^t U^m_\varepsilon$ in \eqref{eqn_help4}, and applying the Cauchy-Schwarz inequality along with the estimate \eqref{eqn_Ubd3}, we obtain
     \begin{align*}
        \|D_+^t U^m_\varepsilon\|^2 &\leq C \|U^m_\varepsilon\|_{1+\alpha}^2 \|D_+^t U^m_\varepsilon\| + \varepsilon^2 \|U^m_\varepsilon\|_{1+\alpha}^2 \|D_+^t U^m_\varepsilon\|,
    \end{align*}
    implies
    \begin{align*}
        \|D_+^t U^m_\varepsilon\| \leq C \|U^m_\varepsilon\|_{1+\alpha}^2 + \varepsilon^2 \|U^m_\varepsilon\|_{1+\alpha}^2 \leq C.
    \end{align*}
    This completes the proof.
\end{proof}

\subsection{Existence, Uniqueness, and Convergence Analysis}
We will prove that the approximations obtained by the scheme \eqref{eqn_fdfsgscheme} converge to the classical solution of the fractional KdV equation \eqref{fkdv} for initial data $u_0\in H^{1+\alpha}_p(I),~\alpha \in [1,2]$. This convergence result can also be interpreted as an existence and uniqueness result for the solution of the fractional KdV equation \eqref{fkdv}.
We first interpolate the approximation $U^m_\varepsilon$ in $[0,\bar T]$ as follows:
\begin{equation}\label{intrpU}
    U_{\varepsilon,N}(x,t) = 
    \begin{cases}
        (1-\theta_m(t))U^{m-\frac{1}{2}}_\varepsilon(x) + \theta_m(t)U^{m+\frac{1}{2}}_\varepsilon(x), & t \in [t_{m-\frac{1}{2}}, t_{m+\frac{1}{2}}), \quad m \geq 1, \\
        \left(1-\frac{t}{\Delta t/2}\right)U^{0}_\varepsilon(x) + \frac{t}{\Delta t/2}U^{\frac{1}{2}}_\varepsilon(x), & t \in (0, t_{\frac{1}{2}}),
    \end{cases}
\end{equation}
where $\theta_m (t) = \frac{1}{\Delta t}(t-t_{m-\frac{1}{2}})$.   {The notation \(U_{\varepsilon,N}(\cdot,t)\in V_N\) for $t\in[0,\bar T]$ denotes the interpolation of approximate solution $\{U_\varepsilon^m\}_{m\Delta t\leq \bar T}$, where \(N\) represents the spatial discretization parameter.}
We show that the approximation $U_{\varepsilon,N}$ converges uniformly to the classical solution of the fractional KdV equation \eqref{fkdv}. More precisely, the next theorem can be seen as the constructive proof of the existence and uniqueness result for the fractional KdV equation \eqref{fkdv}.
\begin{theorem}\label{Con_theorem}
    Let $u_0 \in H_p^{1+\alpha}(I)$ with $\alpha \in [1,2]$ and $T > 0$. Then there exists a unique solution $u \in C([0,T];H_p^{1+\alpha}(I))\cap C^1([0,T];L^2_p(I)) $ to the fractional KdV equation \eqref{fkdv}. 
\end{theorem}

The proof of Theorem \ref{Con_theorem} follows from the lemma and the remark  below, which establishes the existence of the limit of $U_{\varepsilon,N}$ as $N \to \infty$ and asserts that the identified limit is the unique solution of the fractional KdV equation \eqref{fkdv}.

\begin{lemma}\label{conv_lemma}
    Let $U_{\varepsilon,N}$ be an approximate solution of the fractional KdV equation \eqref{fkdv} obtained via the fully discrete scheme \eqref{eqn_fdfsgscheme}, assuming the CFL condition \eqref{eqn_cfl}.  {Then there exist a finite time $\bar T$} such that the following bounds hold
    \begin{align}
        \label{bd1}\norm{U_{\varepsilon,N}(\cdot,t)} &\leq C,\\
        \label{bd2}\norm{U_{\varepsilon,N}(\cdot,t)_t} &\leq C,\\
        \label{bd3}\norm{U_{\varepsilon,N}(\cdot,t)}_{1+\alpha} &\leq C,
    \end{align}
     for all $t\leq\bar T$, where $C = C(\norm{u_0}_{1+\alpha}, \varepsilon, \alpha, \bar T)$. Moreover, the approximation $U_{\varepsilon,N}$ converges to the unique solution of the fractional KdV equation \eqref{fkdv} in $C([0,\bar T];H_p^{1+\alpha}(I))\cap C^1([0,\bar T];L^2_p(I))$.
\end{lemma}

\begin{proof}
   Using the triangle inequality in the definition \eqref{intrpU} of $U_{\varepsilon,N}$, we have
    \begin{equation*}
        \norm{U_{\varepsilon,N}(\cdot,t)} \leq \|U^{m-\frac{1}{2}}_\varepsilon\| + \|U^{m+\frac{1}{2}}_\varepsilon\| \leq C,
    \end{equation*}
    using \eqref{eqn_L2stabfd}, for all $m \geq 0$ and $m\Delta t < \bar T$. The time derivative is given by
    \begin{equation*}
        (U_{\varepsilon,N}(x,t))_t = D_+^tU^{m-\frac{1}{2}}_\varepsilon(x) \qquad t \in [t_{m-\frac{1}{2}}, t_{m+\frac{1}{2}}),
    \end{equation*}
    and thus, for $m \geq 1$ and $m\Delta t <\bar T$,
    \begin{equation}
        \norm{U_{\varepsilon,N}(\cdot,t)_t} = \|D_+^tU^{m-\frac{1}{2}}_\varepsilon\| \leq C, \qquad t \in [t_{m-\frac{1}{2}}, t_{m+\frac{1}{2}}),
    \end{equation}
    using \eqref{eqn_Ubd2}, and similarly for $t \in (0, t_{1/2})$. Lastly, using \eqref{eqn_Ubd3} and the triangle inequality, we obtain
    \begin{equation*}
        \norm{U_{\varepsilon,N}(\cdot,t)}_{1+\alpha} \leq \|U^{m-\frac{1}{2}}_\varepsilon\|_{1+\alpha} + \|U^{m+\frac{1}{2}}_\varepsilon\|_{1+\alpha} \leq C, \qquad t \in [t_{m-\frac{1}{2}}, t_{m+\frac{1}{2}}),
    \end{equation*}
    and similarly for $t \in (0, t_{1/2})$. This completes the proof of \eqref{bd1}-\eqref{bd3}.

    The bound \eqref{bd2} implies that $U_{\varepsilon,N} \in \text{Lip}([0,\bar T];L^2_p(I))$ for all $N$. From \eqref{bd1}, we employ the Arzelà-Ascoli theorem, which ensures the sequential compactness of $\{U_{\varepsilon,N}\}_{N\in\N}$ in $C([0,\bar T];L^2_p(I))$. Consequently, there exists a subsequence $N_k$ such that
    \begin{equation}\label{unifcon}
        U_{\varepsilon,N_k} \to \bar u \text{ uniformly in } C([0,\bar T];L^2_p(I)) \text{ as } N_k \to \infty.
    \end{equation}

    We now show that $\bar u$ is the unique solution of the fractional KdV equation \eqref{fkdv}. First, we claim that $\bar u$ satisfies
    \begin{equation}\label{eqn_wksoln}
       \int_0^{\bar T}\int_I \Big( \bar u  \varphi_t + 3\bar u^2\varphi_x  -\varepsilon^2 \bar u \mathcal{D}^\alpha \varphi_x\Big)\,dx\,dt + \int_I \varphi(x,0)u_0\,dx = 0,
    \end{equation}
    for all  {$\varphi \in C_c^1([0,\bar T]; H_p^{1+\alpha}(I))$.}

For $\varphi(\cdot,t_m) \in H^{1+\alpha}_p(I)$ and $t_m <\bar T$, we choose the test function $\phi = P_N\varphi(\cdot,t_m)$ in \eqref{eqn_fdfsgscheme} and summing it over $m$ after multiplying by $\Delta t$ to obtain
    \begin{align}\label{est_help}
        \nonumber \Delta t \sum\limits_{m\Delta t<\bar T} \left(D_+^tU^{m}_\varepsilon , P_N\varphi(\cdot,t_m)\right) &= -3\Delta t \sum\limits_{m\Delta t<\bar T} \left((U^{m+\frac{1}{2}}_\varepsilon)^2_x , P_N\varphi(\cdot,t_m)\right)\\&\qquad + \varepsilon^2 \Delta t \sum\limits_{m\Delta t<\bar T} \left(\mathcal{D}^{\alpha}(U^{m+\frac{1}{2}}_\varepsilon)_x, P_N\varphi(\cdot,t_m)\right).
    \end{align}
    Using \eqref{eqn_ortho} and the uniform convergence of $U_{\varepsilon,N}$ from \eqref{unifcon}, it follows that
    \begin{align}\label{est1}
        \nonumber \Delta t \sum\limits_{m\Delta t<\bar T} \left(D_+^tU^{m}_\varepsilon , P_N\varphi(\cdot,t_m)\right) &= \Delta t \sum\limits_{m\Delta t<\bar T} \left(D_+^tU^{m}_\varepsilon , \varphi(\cdot,t_m)\right) \\& \nonumber= - \Delta t \sum\limits_{m\Delta t<\bar T} \left(U^{m}_\varepsilon , D_+^t\varphi(\cdot,t_m)\right) - (u_0,\varphi(\cdot,0)) \\
        &\to -\int_0^{\bar T} \int_I \bar u \varphi_t\,dx\,dt - \int_I u_0\varphi(x,0)\,dx,
    \end{align}
    as $N\to \infty$. Similarly, using \eqref{eqn_fracsymm} and \eqref{eqn_ortho}, we have
    \begin{align}\label{est2}
        \nonumber\Delta t \sum\limits_{m\Delta t<\bar T} \left(\mathcal{D}^{\alpha}(U^{m+\frac{1}{2}}_\varepsilon)_x, P_N\varphi(\cdot,t_m)\right) &= \Delta t \sum\limits_{m\Delta t<\bar T} \left(\mathcal{D}^{\alpha}(U^{m+\frac{1}{2}}_\varepsilon)_x, \varphi(\cdot,t_m)\right)\\& = \Delta t \sum\limits_{m\Delta t< \bar T} \left(U^{m+\frac{1}{2}}_\varepsilon, \mathcal{D}^{\alpha}\varphi(\cdot,t_m)_x\right)  \to -\int_0^{\bar T} \int_I \bar u \mathcal{D}^{\alpha}\varphi_x\,dx\,dt,
    \end{align}
    as $N\to\infty$. Also, by \eqref{intrpU},
    \begin{align}\label{est3}
        3\Delta t \sum\limits_{m\Delta t<\bar T} \left((U^{m+\frac{1}{2}}_\varepsilon)^2_x , P_N\varphi(\cdot,t_m)\right) \to -\int_0^{\bar T} \int_I 3\bar u^2\varphi_x\,dx\,dt,
    \end{align}
    as $N\to\infty$. Substituting the results of \eqref{est1}, \eqref{est2}, and \eqref{est3} into \eqref{est_help}, we obtain that \eqref{eqn_wksoln} holds, proving that $\bar u$ is a weak solution of \eqref{fkdv}.  Finally the bounds \eqref{bd1}-\eqref{bd3} implies that $\bar u \in C([0,\bar T];H_p^{1+\alpha}(I))\cap C^1([0,\bar T];L^2_p(I))$ is actually a strong solution which satisfy fractional KdV equation \eqref{fkdv} as an $L^2$-identity.
    To establish the uniqueness, let $\bar u$ and $\bar w$ be two solutions of \eqref{fkdv} in $C([0,\bar T];H_p^{1+\alpha}(I)) \cap C^1([0,\bar T];L^2_p(I))$ with initial conditions $\bar u(\cdot,0) = \bar w(\cdot,0) = u_0$. Define $v = \bar u - \bar w$.  {Then $v$ satisfies the following
    \begin{equation}\label{unq1}
\begin{cases}
      v_t + 6(\bar u\bar u_x - \bar w \bar w_x) -\varepsilon^2 \mathcal{D}^{\alpha} v_x = 0, \qquad &(x,t)\in I\times (0,\bar T],\\
      v(x,0) = 0,  & x \in I,\\
      v(-\pi,t) = v(\pi,t), & \forall t\in [0,\bar T].
\end{cases}
\end{equation}}
Taking inner product of $ v$ and \eqref{unq1}, and applying the identity \eqref{eqn_fracortho}, we obtain
\begin{align*}\label{unq2}
   \nonumber \frac{d}{dt}\norm{v(\cdot,t)}^2 &= -12(\bar u\bar u_x-\bar w\bar w_x,v)= -12(v\bar u_x+\bar wv_x,v) = -12(\bar u_x,v^2) + 6(\bar w_x,v^2)\\&\leq 18C \norm{v(\cdot,t)}^2,
\end{align*}
where we have utilized the Sobolev inequality $\|\bar u_x\|_{\infty}\leq C\|\bar u\|_{1+\alpha}$ and the estimates on solution of the fractional KdV equation \eqref{fkdv}. Applying the Gr\"onwall's inequality, we obtain
\begin{equation*}
    \norm{v(\cdot,t)}^2 \leq e^{18C\bar T} \norm{v(\cdot,0)}^2, \qquad \forall t < \bar T.
\end{equation*}
Since $v(\cdot,0) = 0$, it follows that $v(\cdot,t) = 0$ for all $t < \bar T$. This establishes the uniqueness of the solution. Consequently, we conclude that the sequence of approximations $U_{\varepsilon,N}$ converges uniformly to the unique solution of the fractional KdV equation \eqref{fkdv} in $C([0,\bar T];H_p^{1+\alpha}(I)) \cap C^1([0,\bar T];L^2_p(I))$.
\end{proof}

\begin{remark}\label{remk_extnsn}
 {We now adapt the argument of Sj\"oberg  \cite{sjoberg1970korteweg} to establish global existence for arbitrary time $T>0$. The key observation is that the length of the time interval $[0,\bar{T}]$ of existence depends only on the $H^{1+\alpha}_p$ norm of the initial data $u_0$. Since the exact solution of the fractional KdV equation \eqref{fkdv} preserves this norm, the solution can be extended iteratively.}

 {At time $\bar{T}$, we define the new initial data via projection
\[
U^{\overline{T}}_\varepsilon = P_N u(\cdot, \bar{T}).
\]
 This allows us to extend the solution to the interval $[\bar{T}, 2\bar{T}]$. Repeating this process inductively, we obtain a solution that exists for all $t > 0$, as the $H^{1+\alpha}_p$ norm remains uniformly bounded in time.}
\end{remark}

The above Lemma \ref{conv_lemma} and Remark \ref{remk_extnsn} together demonstrate the existence and uniqueness of the solution to the fractional KdV equation \eqref{fkdv} in \(C([0,T];H_p^{1+\alpha}(I)) \cap C^1([0,T];L^2_p(I))\), completing the constructive proof of Theorem \ref{Con_theorem}.

\subsection{Error Estimate}

To obtain the optimal order of convergence of the fully discrete scheme \eqref{eqn_fdfsgscheme}, we employ the standard approach and define the error function $\mathcal{E}_N^m \in V_N$ by
\begin{equation}
    \mathcal{E}_N^m = P_N u(\cdot,t_m) - U^m_\varepsilon \quad \text{so that} \quad \mathcal{E}_N^0(x) = P_N u_0 - P_N u_0 = 0.
\end{equation}

We will demonstrate that the fully discrete scheme \eqref{eqn_fdfsgscheme} is spectral accurate of order $N^{-r}$ for initial data in Sobolev spaces $H^r_p(I)$ and exponentially accurate for analytic initial data.

\begin{theorem}
    Let $u_0 \in H^r_p(I)$ with $r\geq 3$ and let $u$ be the exact solution of the fractional KdV equation \eqref{fkdv}. Let $U_{\varepsilon,N}$ be the approximate solution obtained by the fully discrete scheme \eqref{eqn_fdfsgscheme} {and assume that the CFL condition \eqref{eqn_cfl} holds.} Then the following estimate holds
    \begin{equation}
        \norm{u(t_m) - U^m_{\varepsilon}} \leq C(N^{-r} + \Delta t^2).
    \end{equation}
\end{theorem}

\begin{proof}
    For all $\phi \in V_N$, we have
\begin{align}\label{err1}
   \nonumber \Bigg(\frac{\mathcal{E}^{m+1}_N - \mathcal{E}^{m}_N}{\Delta t}, \phi\Bigg) &= \Bigg(\frac{P_N u(t_{m+1}) - P_N u(t_m)}{\Delta t}, \phi\Bigg) - \Bigg(\frac{U^{m+1}_\varepsilon - U^m_\varepsilon}{\Delta t}, \phi\Bigg)\\
   \nonumber
    &= (\xi^{m+\frac{1}{2}}, \phi) +  (u_t(t_{m+\frac{1}{2}}), \phi) + 3 \Big((U^{m+\frac{1}{2}}_\varepsilon)_x^2, \phi\Big) - \varepsilon^2 \Big(\mathcal{D}^\alpha (U^{m+\frac{1}{2}}_\varepsilon)_x, \phi\Big)\\
    & = (\xi^{m+\frac{1}{2}}, \phi) + 3 \Big((U^{m+\frac{1}{2}}_\varepsilon)_x^2 - u(t_{m+\frac{1}{2}})^2_x, \phi\Big) - \varepsilon^2 \Big(\mathcal{D}^\alpha (U^{m+\frac{1}{2}}_\varepsilon - u(t_{m+\frac{1}{2}}))_x, \phi\Big),
\end{align}
    where $$\xi^{m+\frac{1}{2}} = \frac{u(t_{m+1}) - u(t_m)}{\Delta t} - u_t(t_{m+\frac{1}{2}}).$$
Let us define the error at the average time $t_{m+\frac{1}{2}}$ by
    \begin{align*}
        U^{m+\frac{1}{2}}_\varepsilon - P_N u(t_{m+\frac{1}{2}}) =  -\mathcal{E}_N^{m+\frac{1}{2}} + P_N \theta^{m+\frac{1}{2}},
    \end{align*}
    where $$\mathcal{E}_N^{m+\frac{1}{2}} = \frac{\mathcal{E}_N^{m+1} + \mathcal{E}_N^{m}}{2},\quad \text{and} \quad \theta^{m+\frac{1}{2}} = \frac{u(t_{m+1}) + u(t_m)}{2} - u(t_{m+\frac{1}{2}}).$$
Using the Taylor's expansion, it is straightforward to show that 
\begin{equation}\label{taylors_est}
    \|P_N \theta^{m+\frac{1}{2}}\|^2_l \leq C \Delta t^3 \sup\limits_{t<T} \|u_{tt}(t)\|^2_l, \quad \text{and} \quad \|\xi^{m+\frac{1}{2}}\|^2 \leq C \Delta t^3 \sup\limits_{t<T} \|u_{ttt}(t)\|^2,
\end{equation}
for all $l \geq r-2$. By taking $\phi = \mathcal{E}^{m+\frac{1}{2}}_N$ in \eqref{err1}, we obtain
\begin{align}
    \nonumber \norm{\mathcal{E}^{m+1}_N}^2 &\leq \norm{\mathcal{E}^{m}_N}^2 + 2\Delta t\Big[C \Delta t^3 + \frac{1}{8} \|\mathcal{E}^{m+\frac{1}{2}}_N\|^2 + 3 \Big((U^{m+\frac{1}{2}}_\varepsilon)_x^2 - (u(t_{m+\frac{1}{2}})^2)_x, \mathcal{E}^{m+\frac{1}{2}}_N\Big)\\&\qquad \label{erreqn} + \varepsilon^2 C \Delta t^3 + \frac{1}{4} \|\mathcal{E}^{m+\frac{1}{2}}_N\|^2\Big].
\end{align}
We estimate the nonlinear term in the above equation as follows
\begin{align}
    \nonumber 3\Big(&(U^{m+\frac{1}{2}}_\varepsilon)_x^2 - (u(t_{m+\frac{1}{2}})^2)_x, \mathcal{E}^{m+\frac{1}{2}}_N\Big)\\& \nonumber= 6\Big(U^{m+\frac{1}{2}}_\varepsilon (U^{m+\frac{1}{2}}_\varepsilon)_x - P_N (u(t_{m+\frac{1}{2}}) u(t_{m+\frac{1}{2}})_x), \mathcal{E}^{m+\frac{1}{2}}_N\Big)\\& \nonumber= 6\Big( (U^{m+\frac{1}{2}}_\varepsilon - P_N u(t_{m+\frac{1}{2}}))(U^{m+\frac{1}{2}}_\varepsilon - P_N u(t_{m+\frac{1}{2}}))_x, \mathcal{E}^{m+\frac{1}{2}}_N\Big) + 6\Big(P_N u(t_{m+\frac{1}{2}})(U^{m+\frac{1}{2}}_\varepsilon)_x, \mathcal{E}^{m+\frac{1}{2}}_N\Big)\\
    & \nonumber\quad- 6\Big(P_N u(t_{m+\frac{1}{2}})P_N u(t_{m+\frac{1}{2}})_x, \mathcal{E}^{m+\frac{1}{2}}_N\Big) + 6\Big(P_N u_x(t_{m+\frac{1}{2}})(U^{m+\frac{1}{2}}_\varepsilon - P_N u(t_{m+\frac{1}{2}})), \mathcal{E}^{m+\frac{1}{2}}_N\Big)\\
    & \nonumber
    \leq 12C \Delta t^3 + \frac{3}{4} \|\mathcal{E}^{m+\frac{1}{2}}_N\|^2 + 12 \|u_x(t_{m+\frac{1}{2}})\|_\infty \|\mathcal{E}^{m+\frac{1}{2}}_N\|^2 + 12 \|u_x(t_{m+\frac{1}{2}})\|_\infty \Big(\|\mathcal{E}^{m+\frac{1}{2}}_N\|^2 + C \Delta t^3\Big)\\
     \label{help2err} &\leq 12C \Delta t^3 + C_0 \|\mathcal{E}^{m+\frac{1}{2}}_N\|^2,
\end{align}
  {where $C$ and $C_0$ are constants independent of $\Delta t$, and the result follows from the use of \eqref{taylors_est}, Lemma \ref{lemma_pbound}, and \eqref{eqn_cfl}.}
 Substituting \eqref{help2err} into \eqref{erreqn}, we obtain
     \begin{align*}
    \norm{\mathcal{E}^{m+1}_N}^2 \leq \norm{\mathcal{E}^{m}_N}^2 + 2\Delta t\Big[C\Delta t^3 +C_1\|\mathcal{E}^{m+\frac{1}{2}}_N\|^2\Big] \leq (1+C_1\Delta t) \norm{\mathcal{E}^{m}_N}^2 + C\Delta t^4 + C_1\Delta t\|\mathcal{E}^{m+1}_N\|^2,
   \end{align*}
where constant $C$ and $C_1$ are generic and independent of $\Delta t$. We choose $\Delta t$ sufficiently small such that $1-\Delta t C_1 \geq \frac{1}{2}$ and use the fact that $\mathcal{E}^0_N = 0$. The above estimate implies
\begin{align*}
   \norm{\mathcal{E}^{m+1}_N}^2 &\leq \frac{1+C_1\Delta t}{1-C_1\Delta t} \norm{\mathcal{E}^{m}_N}^2 +  \frac{C\Delta t^4}{1-C_1\Delta t}
    \leq \left(\frac{1+C_1\Delta t}{1-C_1\Delta t}\right)^{m+1}\norm{\mathcal{E}^{0}_N}^2 + e^{CT}\Delta t^4 
    \leq e^{CT}\Delta t^4.
\end{align*}
Furthermore, since $u_0\in H_p^r(I)$, we have
\begin{equation}\label{errest2}
    \norm{U^m_\varepsilon - P_Nu(t_{m})} \leq \norm{P_Nu(t_m)-u(t_m)} + \norm{\mathcal{E}^m_N}\leq C(N^{-r} + \Delta t^2),
\end{equation}
where $C$ is dependent on $\alpha, \varepsilon, T,$ and $\norm{u_0}_r$. This completes the proof.
\end{proof}

We observe that, if the the initial data $u_0$ is analytic and hence the solution of fractional KdV equation \eqref{fkdv} is analytic, then the estimate \eqref{errest2} and the Proposition \ref{Prop_exp} implies
\begin{equation}\label{errest3}
    \norm{U^m_\varepsilon - P_Nu(t_{m})} \leq \norm{P_Nu(t_m)-u(t_m)} + \norm{\mathcal{E}^m_N}\leq C(e^{-cN} + \Delta t^2).
\end{equation}

\section{Zero Dispersive Limit}\label{sec4}
In this section, we investigate the behavior of solutions to the fractional KdV equation \eqref{fkdv} as the dispersion parameter $\varepsilon$ tends to zero. This regime is known as the zero dispersive limit. The analysis is divided into two cases: for time $ t < t_c $, where $ t_c $ is the critical times corresponding to the onset of gradient catastrophe in the inviscid Burgers' equation, and for times $ t > t_c $.

For $ t < t_c $, the limiting equation
\begin{equation}\label{eqn_burgers}
    u_t + 6u u_x = 0,
\end{equation}
with smooth initial data $ u_0 $, admits a smooth solution. Previously, we established that for each $ \varepsilon > 0 $, $ T > 0 $, and initial data $ u_0 \in H^{1+\alpha}_p(I) $, there exists a unique solution $ u(x,t;\varepsilon) $ to the fractional KdV equation \eqref{fkdv} such that
\begin{equation}
    u(x,t;\varepsilon) \in C([0,T];H^{1+\alpha}_p(I)) \cap C^1([0,T];L^2_p(I)).
\end{equation}
Our goal is to examine the limiting behavior of $ u(x,t;\varepsilon) $ as $ \varepsilon \to 0 $, particularly the quantity $ \lim\limits_{\varepsilon \to 0}\lim\limits_{N \to \infty} U_{\varepsilon,N} = \lim\limits_{\varepsilon \to 0} u(x,t;\varepsilon) $.

\subsubsection*{Case 1: $ t < t_c $}

For $ t < t_c $, the solution remains smooth, and the convergence of $ u(x,t;\varepsilon) $ to the solution $ v^b(x,t) $ of the limiting equation \eqref{eqn_burgers} is straightforward. This is captured in the following theorem.

\begin{theorem}\label{thm:main-cont}
    Let $ T_c $ be any positive time less than the critical time $ t_c $. Let $ u(x,t;\varepsilon) \in C([0,T_c];H^{1+\alpha}_p(I)) \cap C^1([0,T_c];L^2_p(I)) $ be the unique solution of the fractional KdV equation \eqref{fkdv}. Then the map $ \varepsilon \mapsto u(x,t;\varepsilon) $ is continuous from $ \mathbb{R}^+ $ to $ C([0,T_c];H^{1+\alpha}_p(I)) \cap C^1([0,T_c];L^2_p(I)) $. Moreover, we have
    \begin{equation}\label{limitofu}
        \lim\limits_{\varepsilon \to 0} \norm{u(t;\varepsilon) - v^b(t)}_{1+\alpha} = 0,
    \end{equation}
    uniformly in $ t \in [0,T_c] $, where $ v^b(x,t) $ satisfies the limiting equation
    \begin{equation}\label{eqn_burger}
        v^b_t + 6v^b v^b_x = 0, \qquad v^b(x,0) = u_0(x).
    \end{equation}
\end{theorem}
To prove this, we first establish an abstract continuity result for parameter-dependent quasilinear equations. Borrowing the idea from \cite[Theorem 4-5]{masoero2013semiclassical} and \cite[Theorem 7]{kato2006quasi}, we have the following result.

\begin{theorem}[Abstract parameter continuity]\label{thm:abstract-cont}
 {Let \(X,Y\) be Hilbert spaces with \(Y \subset X\) densely and continuously embedded. Consider the family of quasilinear problems
\begin{equation}\label{abstract-eq}
\frac{du}{dt} = A(u;\varepsilon)u, \quad u(0) = u_0 \in Y, \quad \varepsilon \in \mathbb{R}^+,
\end{equation}
where \(A: Y \times \mathbb{R}^+ \to \mathcal{L}(Y,X)\) and assume that $\|A(u(\varepsilon);\varepsilon)\|_{Y, X}$ is uniformly bounded in $u$ and $\varepsilon$. Furthermore, assume that}

\begin{enumerate}[label=(H\arabic*)]
    \item  {\textbf{Well-posedness}: For each \(\varepsilon > 0\) and \(T>0\), \eqref{abstract-eq} has a unique solution \(u(\cdot;\varepsilon) \in C([0,T];Y) \cap C^1([0,T];X)\).}
    
    \item  {\textbf{Lipschitz continuity in \(u\)}: There exists \(L > 0\) such that for all \(\varepsilon > 0\),
    \begin{equation}
    \|A(u;\varepsilon) - A(v;\varepsilon)\|_{Y,X} \leq L\|u - v\|_Y.
    \end{equation}}
    
    \item  {\textbf{Parameter regularity}:
    \begin{equation}
    \|A(u;\varepsilon) - A(u;\varepsilon')\|_{Y,X} \leq C|\varepsilon - \varepsilon'|.
    \end{equation}}
    
    \item  {\textbf{Commutator structure}: There exists an isometric isomorphism \(\Lambda: Y \to X\) such that
    \begin{equation}
    \|[\Lambda, A(u;\varepsilon)]\Lambda^{-1}\|_{X} \leq C\|u\|_Y,
    \end{equation}
     where the commutator operator is defined as $[\Lambda, A(u;\varepsilon)]u:= \Lambda A(u;\varepsilon) u-A(u;\varepsilon)\Lambda u$.} 
      (The above relation should be satisfied in the strict sense, including the domain relation. Thus $x\in X$ is in domain of $A(u;\varepsilon)$ if and only if $\Lambda^-x$ is in domain of $A(u;\varepsilon)$ with $A(u;\varepsilon)\Lambda^-x\in Y$. For more details we refer to \cite[Assumption (A2)]{kato2006quasi}.)
\end{enumerate}
 {Then the solution map \(\varepsilon \mapsto u(\cdot;\varepsilon)\) is continuous from \(\mathbb{R}^+\) to \(C([0,T];Y) \cap C^1([0,T];X)\).}
\end{theorem}

\begin{proof}[Proof of Theorem \ref{thm:abstract-cont}]
  {Condition (H1) provides the well-posedness framework, while (H2)-(H4) give the required continuity properties.
By (H1), for each \(\varepsilon > 0\), the solution \(u(\cdot;\varepsilon)\) exists on \([0,T]\). Since \(u_0 \in Y\) is fixed and \(A(\cdot;\varepsilon)\) is locally Lipschitz (from (H2)), a standard energy estimate (using Gr\"onwall's inequality) shows that \(\|u(t;\varepsilon)\|_Y\) is uniformly bounded for all \(t \in [0,T]\) and \(\varepsilon>0\).} 

 {Let \(u(\varepsilon)\) and \(u(\varepsilon')\) be solutions corresponding to parameters \(\varepsilon, \varepsilon' \in \R^+\). The difference \(w = u(\varepsilon) - u(\varepsilon')\) satisfies
\begin{equation}\label{diff-eq}
w_t = A(u(\varepsilon);\varepsilon)w + E,
\end{equation}
where \(E = (A(u(\varepsilon);\varepsilon) - A(u(\varepsilon');\varepsilon'))u(\varepsilon')\).
To estimate \(\|w(t)\|_Y\), apply the isometric isomorphism \(\Lambda: Y \to X\) and compute
\begin{align*}
\frac{1}{2}\frac{d}{dt}\| w\|_Y^2 &= \langle \Lambda w, \Lambda A(u(\varepsilon);\varepsilon) \Lambda^-(\Lambda w) \rangle_X + \langle w,  E \rangle_Y \\
&= \langle \Lambda w, A(u(\varepsilon);\varepsilon)\Lambda w \rangle_X + \langle \Lambda w, [\Lambda,A(u(\varepsilon);\varepsilon)]\Lambda^{-}(\Lambda w) \rangle_X + \langle  w,  E \rangle_Y,
\end{align*} 
where $[\Lambda, A(u;\varepsilon)]\Lambda^{-1} =\Lambda A(u;\varepsilon)\Lambda^{-1}-A(u;\varepsilon)$. Using (H4), the commutator term is bounded as
\begin{align*}
|\langle \Lambda w, [\Lambda,A(u(\varepsilon);\varepsilon)]w \rangle_X| &\leq \|\Lambda w\|_X \|[\Lambda,A(u(\varepsilon);\varepsilon)]\Lambda^{-1}\|_{X} \|\Lambda w\|_X \\
&\leq C\|u(\varepsilon)\|_Y \|w\|_Y^2.
\end{align*}
Furthermore, \(E\) can be decomposed as
\[
E = (A(u(\varepsilon);\varepsilon) - A(u(\varepsilon');\varepsilon))u(\varepsilon') + (A(u(\varepsilon');\varepsilon) - A(u(\varepsilon');\varepsilon'))u(\varepsilon').
\]
By (H2) and (H3), we have:
\[
\|E\|_Y \leq L\|w\|_Y \|u(\varepsilon')\|_Y + C|\varepsilon - \varepsilon'|\|u(\varepsilon')\|_Y.
\]
Note that \(\|u(\varepsilon')\|_Y< C\) by (H1) for each $\varepsilon>0$, so
\[
\|E\|_Y \leq C\|w\|_Y + C|\varepsilon - \varepsilon'|.
\]
Since $\|A(u(\varepsilon);\varepsilon)\|_{Y, X}$ is uniformly bounded in $u$ and $\varepsilon$, from isometry, we have
 $$|\langle \Lambda w, A(u(\varepsilon);\varepsilon)\Lambda w \rangle_X| \leq \|A(u(\varepsilon);\varepsilon)\|_{Y, X}\|\Lambda w\|_X^2 \leq C\|u(\varepsilon)\|_Y\|w\|_Y^2.$$
Combining the above estimates and using the isometry of \(\Lambda\), we end up with
\[
\frac{d}{dt}\|w\|_Y^2 \leq C\|w\|_Y^2 + C|\varepsilon - \varepsilon'|\|w\|_Y.
\]
Using the Young’s inequality we obtain
\[
\frac{d}{dt}\|w\|_Y^2 \leq C\|w\|_Y^2 + C|\varepsilon - \varepsilon'|^2.
\]
The Gr\"onwall's inequality further yields
\[
\|w(t)\|_Y \leq C|\varepsilon - \varepsilon'|e^{CT} \quad \text{for all } t \in [0,T].
\]
This establishes the continuity of \(\varepsilon \mapsto u(\cdot;\varepsilon)\) in \(C([0,T];Y) \cap C^1([0,T];X)\).}
\end{proof}

\begin{proof}[Proof of Theorem \ref{thm:main-cont}]
 {We employ the Kato's theory \cite{kato2006quasi} and rewrite the fractional KdV equation \eqref{fkdv} as 
    \begin{align}\label{A_1}
        &\frac{du}{dt} = A(u;\varepsilon)u, \qquad u(0) = u_0 \in H^{1+\alpha}_p(I), \\
        &A(u;\varepsilon) = -6u\partial_x + \varepsilon^2 \mathcal{D}^\alpha \partial_x,\label{A_2}
    \end{align}
    where $ A(u;\varepsilon) $ can be considered as the perturbed operator of the Burgers' operator corresponding to $ \varepsilon > 0 $.
We verify the hypotheses of Theorem \ref{thm:abstract-cont} for the fractional KdV operator $A(u;\varepsilon)$ for $u\in H_p^{1+\alpha}$ with $X = L^2_p(I)$ and $Y = H^{1+\alpha}_p(I)$, and $X$ and $Y$ are Hilbert Space. The well-posedness of \eqref{A_1} was established in Section \ref{sec3} which verifies the hypothesis (H1) of Theorem \ref{thm:abstract-cont}. The fractional operator  $\mathcal{D}^\alpha\partial_x$ is a bounded operator from  $H^{1+\alpha}_p(I)$ to $L^2_p(I)$. For $u\in H^{1+\alpha}_p(I)$, the nonlinear term $6u\partial_x : H^{1+\alpha}_p(I) \to H^\alpha_p(I) \hookrightarrow L^2_p(I)$ is bounded by
\begin{equation}
\|u\partial_x v\| \leq \|u\|_{\infty}\|\partial_x v\| \leq C\|u\|_{1+\alpha}\|v\|_{1+\alpha}
\end{equation}
using the Sobolev embedding 
$\|u\|_{\infty}\leq C \|u\|_{1+\alpha}$. Hence using these boundedness of operators, for $u,v \in H^{1+\alpha}_p(I)$
 \begin{align}\label{est110}
 \|A(u;\varepsilon) - A(v;\varepsilon)\|_{H_p^{1+\alpha},L_p^2} \leq C\|u-v\|_{1+\alpha}.
 \end{align}
 This implies that $A(u;\varepsilon)$ is a Lipschitz operator from $H_p^{1+\alpha}(I)$ to $L^2_p(I)$, for $\varepsilon\in \R$. For fixed $u \in H^{1+\alpha}_p(I)$,
    \begin{equation}\label{est10}
    \|A(u;\varepsilon) - A(u;\varepsilon')\|_{H_p^{1+\alpha},L_p^2} \leq C|\varepsilon^2 - \varepsilon'^2|.
    \end{equation}
 Let $\Lambda^{1+\alpha} = (1-\partial_x^2)^{\frac{1+\alpha}{2}}$ be the isometric isomorphism $H^{1+\alpha}_p(I) \to L^2_p(I)$. Then
\begin{equation}\label{est11}
B(u;\varepsilon) := (\Lambda^{1+\alpha}A(u;\varepsilon) - A(u;\varepsilon)\Lambda^{1+\alpha})\Lambda^{-(1+\alpha)}= [\Lambda^{1+\alpha},A(u;\varepsilon)]\Lambda^{-(1+\alpha)}
\end{equation}
where $ B(u;\varepsilon)$ satisfies the commutator estimate (\cite[Lemma 1]{masoero2013semiclassical}, \cite[Lemma A.2]{kato2006quasi})
\begin{equation}\label{est13}
\|B(u;\varepsilon)\| \leq C(\alpha)\|u\|_{1+\alpha}.
\end{equation}
  Finally, estimates \eqref{est110}-\eqref{est13} verify conditions (H2)-(H4) of Theorem \ref{thm:abstract-cont}, and hence the map $\varepsilon \mapsto u(x,t;\varepsilon)$ is continuous from $\R^+$ to $C([0,T_c];H^{1+\alpha}_p(I)) \cap C^1([0,T_c];L^2_p(I))$.}
 
  {Afterwards, we consider the solution operator $\Phi_\varepsilon$ for the fractional KdV equation, which maps initial data to solutions at time $t$
\[
\Phi_\varepsilon(t) u_0(x) = u(x,t;\varepsilon).
\]
We have that $\Phi_\varepsilon(t)$ depends continuously on $\varepsilon$ in the operator norm topology of $\mathcal{L}(H_p^{1+\alpha}(I), H^{1+\alpha}_p(I))$, uniformly for $t \in [0,T_c]$. The key observation at $\varepsilon = 0$ is that the fractional KdV operator reduces to $A(u;0) = -6u\partial_x$, which generates exactly the Burgers' equation flux. Thus, at $\varepsilon = 0$, the operator $\Phi_0(t)$ exactly corresponds to the solution map for the Burgers' equation \eqref{eqn_burger}
\[
\Phi_0(t)u_0(x) = v^b(x,t).
\]
The classical theory of Burgers' equation guarantees that $v^b \in C([0,T_c]; H^{1+\alpha}_p(I))$ for initial data in $H^{1+\alpha}_p(I)$ when $T_c < t_c$. Finally, the uniform convergence follows from the operator norm continuity
\[
\sup_{t \in [0,T_c]} \|\Phi_\varepsilon(t)u_0 - \Phi_0(t)u_0\|_{1+\alpha} \leq \sup_{t \in [0,T_c]} \|\Phi_\varepsilon(t) - \Phi_0(t)\|_{H_p^{1+\alpha},H_p^{1+\alpha}} \|u_0\|_{1+\alpha}.
\]
The right-hand side tends to zero as $\varepsilon \to 0$ by the uniform continuity of $\varepsilon \mapsto \Phi_\varepsilon$. This proves the uniform convergence \eqref{limitofu}.}

\end{proof}

 \begin{figure}
    \centering
    \includegraphics[width=0.9 \linewidth, height=9cm]{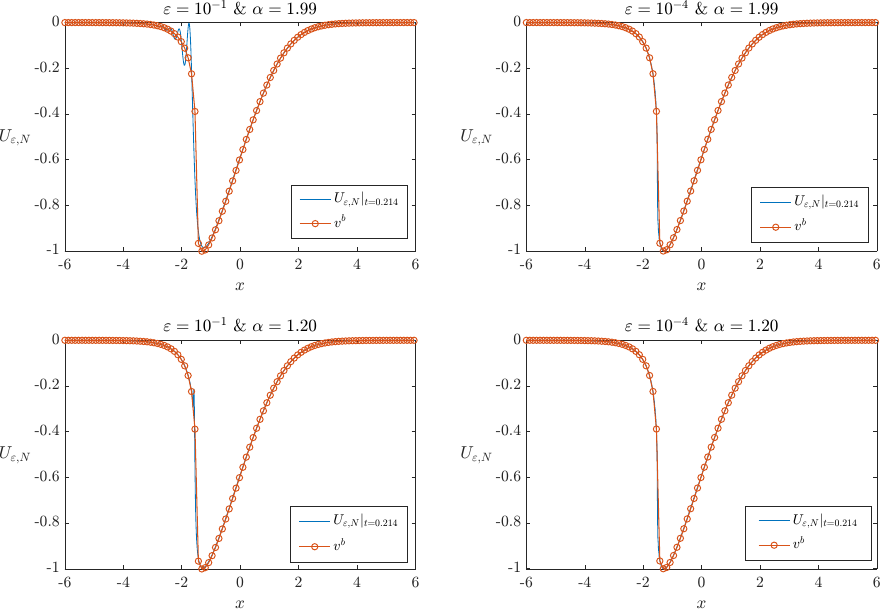}
    \caption{ The numerical approximation $U_{\varepsilon,N}$  {with $N=2^{16}$} of the fractional KdV equation \eqref{fkdv} with $\alpha=1.999$ and $\alpha=1.20$ at time $t = 0.214<t_c$ for different dispersive coefficients $\varepsilon$ and the solution $v^b$ of the limiting equation \eqref{eqn_burgers}.}
    \label{fig:breaktime}
\end{figure}

\subsubsection*{Case 2: $ t > t_c $} 
After the gradient catastrophe arises in the limiting equation, the behavior of the solution to the fractional KdV equation \eqref{fkdv} becomes more intricate, as it does not converge to the solution of the limiting equation \eqref{eqn_burger} when $\varepsilon \to 0$. We examine the limit of the fractional KdV solution as $\varepsilon \to 0$ for times $ t > t_c $. A numerical comparison is made between the approximations $ U_{\varepsilon,N} $ for $ \alpha $ close to $2$ and the corresponding KdV equation. 

Building on the foundational work of Lax and Levermore \cite{lax1983small,lax1983small2,lax1983small3}, Venakides \cite{venakides1987zero}, and Grava and Klein \cite{grava2012numerical, grava2007numerical, klein2015numerical}, it is known that after the time $ t_c $, oscillatory zone $[x^-(t), x^+(t)]$ arises, which is independent of $\varepsilon$ for single hump, rapidly decreasing initial data. The boundaries $ x^-(t) $ and $ x^+(t) $ are determined by the initial data and satisfy $ x^-(t_c) = x^+(t_c) = x_c $, where $ x_c $ denotes the spatial coordinate of the gradient catastrophe at the break time $ t_c $ of the limiting equation. Prior to implementing the numerical scheme \eqref{eqn_fdfsgscheme} for solving the zero dispersive limit of the limiting equation for $ t > t_c $, we introduce the notion of asymptotic solutions for the fractional KdV equation \eqref{fkdv} across different values of $\alpha$. This introduction allows us to evaluate the performance of the devised scheme \eqref{eqn_fdfsgscheme} by comparing it with the asymptotic behavior of the solution of the fractional KdV equation \eqref{fkdv} as $\varepsilon$ approaches zero.

Outside the oscillatory zone $[x^-(t), x^+(t)]$, the asymptotic solution $ u(x,t;\varepsilon) $ as $\varepsilon \to 0$ for all $\alpha \in [1,2]$ converges to the solution of the limiting equation. Within the oscillatory zone, for $\alpha = 2$ and sufficiently small $\varepsilon$, the solution $ u(x,t;\varepsilon) $ is approximately described by the elliptic solution of the KdV equation
\begin{equation}\label{kdvellson}
    u_t + 6uu_x + \varepsilon^2 u_{xxx} = 0,
\end{equation}
which is given by \cite{lax1983small, lax1983small2, lax1983small3, venakides1987zero, grava2007numerical} as:
\begin{equation}\label{asyptoticsol}
    u(x,t;\varepsilon) \approx \Tilde{u} + 2\varepsilon^2 \frac{\partial^2}{\partial x^2} \log \mu\left(\frac{\sqrt{\beta_1 - \beta_3}}{2\varepsilon K(s)} \left[x - 2t(\beta_1 + \beta_2 + \beta_3) - q\right]; \mathcal{T}\right),
\end{equation}
where $\Tilde{u} = \Tilde{u}(x,t)$ is the weak limit of $u(x,t;\varepsilon)$ as $\varepsilon \to 0$ \cite{lax1983small}, has the following form
\begin{equation*}
    \Tilde{u}= \beta_1+\beta_2+\beta_3+2\beta, \text{ and } \beta = -\beta_1 + (\beta_1-\beta_3)\frac{E(s)}{K(s)}, 
\end{equation*}
with $K(s)$ and $E(s)$ being the first and second kind complete elliptic integrals. Note that 
\begin{equation*}
    K'(s) = K(\sqrt{1-s^2}), \qquad \mathcal{T} = i \frac{K'(s)}{K(s)}, \qquad s^2 = \frac{\beta_3-\beta_1}{\beta_3-\beta_1}.
\end{equation*}
Furthermore, the Jacobi elliptic function $\mu$ is defined by the Fourier series
\begin{equation*}
    \mu(\xi;\mathcal{T}) = \sum\limits_{n\in\Z}e^{\pi i n^2\mathcal{T}+2\pi i n\xi}.
\end{equation*}
The quantities $\beta_i(x,t)$, $i=1,2,3$, evolve according to the Whitham equations \cite{whitham2011linear}
\begin{equation}\label{whitham}
    \frac{\partial}{\partial t}\beta_i + v_i\frac{\partial}{\partial x}\beta_i = 0,\quad i=1,2,3, \quad\text{and}\quad v_i = 4 \frac{\prod_{j\neq i}(\beta_i-\beta_k)}{\beta_i+\beta} + 2 (\beta_1+\beta_2+\beta_3), \quad i=1,2,3.
\end{equation}
The formula for $q$ in \eqref{kdvellson} is given by \cite{grava2007numerical}
\begin{align*}
    q(\beta_1,\beta_2,\beta_3) = \frac{1}{2\sqrt{2}\pi}\int_{-1}^1\int_{-1}^1\frac{f_-(A)}{\sqrt{1-\theta}\sqrt{1-\gamma^2}}\,d\theta\,d\gamma,
\end{align*}
where $f_-$ is the inverse function in the decreasing part of the initial data and $A = \frac{1+\theta}{2}(\frac{1+\gamma}{2}\beta_1+\frac{1-\gamma}{2}\beta_2)+\frac{1-\theta}{2}\beta_3$.

Now we consider an explicit example to illustrate the numerical scheme \eqref{eqn_fdfsgscheme} and asymptotic solution \eqref{kdvellson}. 
\subsection*{Example 4.1}
    For $\alpha \in [1,2]$, consider the fractional KdV equation \eqref{fkdv} with initial condition
    \begin{equation}\label{usualkdvinitial}
       u_0(x) = -\sech^2(x),
    \end{equation}
    for $(x,t) \in \R \times (0,T)$. The gradient catastrophe point $(x_c, t_c, u_c)$, where $u_c = u(x_c, t_c)$, is analytically given by:
    \begin{equation*}
        t_c = \frac{1}{\displaystyle\max_{x\in\R}[-6u_0'(x)]}= \frac{\sqrt{3}}{8}, \quad x_c = -\frac{\sqrt{3}}{2} + \log\left(\frac{\sqrt{3}-1}{\sqrt{2}}\right), \quad u_c = -\frac{2}{3}.
    \end{equation*}
    The asymptotic solution of the usual KdV equation \eqref{kdvellson} in the oscillatory zone $[x^-(t), x^+(t)]$ is given by \eqref{kdvellson}. Outside this zone, the solution is described by:
    \begin{equation}\label{hopfsoln}
        u(x,t) = u_0(\xi), \qquad x = 6tu_0(\xi) + \xi.
    \end{equation}
    Numerical simulations of the scheme \eqref{eqn_fdfsgscheme} are performed at times $t = 0.2 < t_c$ and $t = 0.4 > t_c$ to estimate the error for different values of $\varepsilon$. The error $E(\varepsilon)$ for fixed time $t$ is defined as:
    \begin{equation*}
        E(\varepsilon) :=\sup\limits_{x\in[-L,L]}|u(x,t;\varepsilon)-U_{\varepsilon,N}(x,t)|.
    \end{equation*}
    A time step \(\Delta t = 1/(8N\|u_0\|_{\infty})\) is used, with period \(L=6\) and \(N = 2^{16}\) to ensure sufficient resolution.  {Figure \ref{fig:breaktime} presents the result of the numerical simulation for \(t < t_c\), where the approximate solution is compared with the exact solution of \eqref{eqn_burgers} associated with the initial data \eqref{usualkdvinitial}.}  Figures \ref{fig:errzdl}–\ref{fig:alphavar} present the results of these numerical simulations for \(t > t_c\). Table \ref{tab:error_tableeps} represents the error $E(\varepsilon)$ with respect to small $\varepsilon$ for time $t<t_c$ and $t>t_c$.
    We numerically verify that for \(\alpha = 1.999 \approx 2\) and small \(\varepsilon^2\), the approximate solution obtained using the scheme \eqref{eqn_fdfsgscheme} closely matches the asymptotic solutions \eqref{asyptoticsol} and \eqref{hopfsoln} of the usual KdV equation \eqref{kdvellson}.

\begin{figure}
    \centering
    \includegraphics[width=0.7 \linewidth, height=7cm]{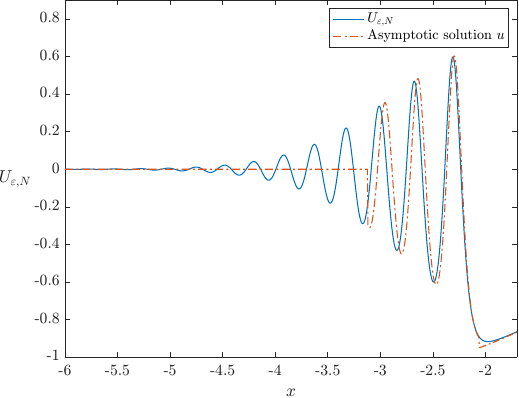}
    \caption{ The numerical approximation $U_{\varepsilon,N}$ of the fractional KdV equation \eqref{fkdv} and asymptotic solution \eqref{hopfsoln} and \eqref{asyptoticsol} with $\alpha=1.999$ {at time $t=0.4$} for the dispersive coefficient $\varepsilon = 10^{-1}$.}
    \label{fig:errzdl}
\end{figure}
\begin{table}[htbp]
        \centering
        \begin{tabular}{||c|c|c|c||}
        \hline
              \multicolumn{2}{||c|}{$t=0.2<t_c$} & \multicolumn{2}{c||}{$t=0.4>t_c$} \\
            \hline
          $\varepsilon$   & $E(\varepsilon)$ & $\varepsilon$  & $E(\varepsilon)$    \\
            \hline
            \hline
              $10^{-1.0}$  &   2.25e-01 & $10^{-1.0}~$  &  7.03e-01 \\  
             $10^{-2.0}$ & 3.41e-02 & $10^{-1.5}~$  & 3.05e-01   \\
              $10^{-2.5}$ & 3.91e-03  &$10^{-2.0}~$  & 8.86e-02  \\
             $10^{-3.0}$&   4.02e-04 &$10^{-2.5}~$  &  2.59e-02 \\  
             $10^{-3.5}$&   4.09e-05 &$10^{-2.8}~$  & 5.01e-03 \\  
             $10^{-4.0}$ &  4.10e-06 &$10^{-3.0}~$  & 1.41e-03 \\
            \hline
        \end{tabular}
        \caption{Error $E(\varepsilon)$ for the scheme \eqref{eqn_fdfsgscheme} at times $t=0.2<t_c$ and  $t=0.4>t_c$ taking $N=2^{16}$.}
        \label{tab:error_tableeps}
 \end{table}
\begin{figure}
    \centering
    \includegraphics[width=0.8 \linewidth, height=7cm]{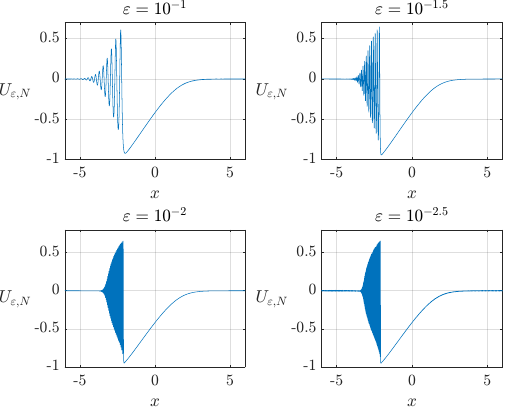}
    \caption{ The numerical approximation $U_{\varepsilon,N}$ of the fractional KdV equation \eqref{fkdv} with $\alpha =1.999$ at the time $t = 0.4$ and for different dispersive coefficients $\varepsilon$.}
    \label{fig:epsvar}
\end{figure}

 \begin{figure}
    \centering
    \includegraphics[width=0.8 \linewidth, height=7cm]{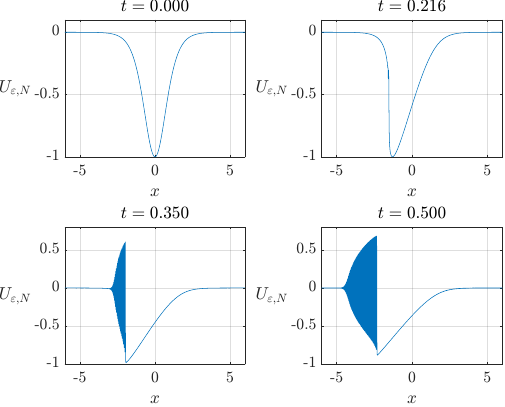}
    \caption{ The numerical approximation $U_{\varepsilon,N}$ of the fractional KdV equation \eqref{fkdv} with $\alpha=1.999$ at different times $t$ for the dispersive coefficient $\varepsilon = 10^{-2}$.}
    \label{fig:timevar}
\end{figure}
\begin{figure}
    \centering
    \includegraphics[width=0.8 \linewidth, height=7cm]{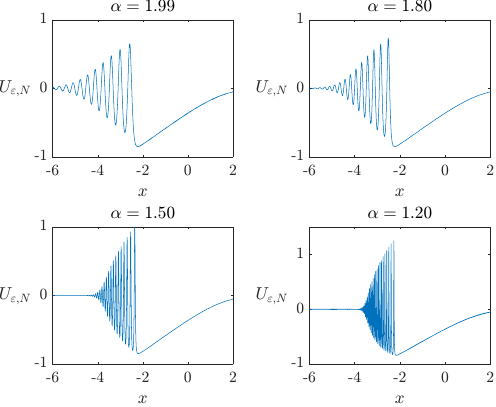}
    \caption{ The numerical approximation $U_{\varepsilon,N}$ of the fractional KdV equation \eqref{fkdv} with the initial data \eqref{usualkdvinitial} at the time $t = 0.5$, and with the coefficient $\varepsilon=10^{-1}$ and for different values of $\alpha$.}
    \label{fig:alphavar}
\end{figure}

\subsection*{Example 4.2}
    We consider the fractional KdV equation 
    \begin{equation}\label{usualkdvbo}
        u_t + uu_x - \varepsilon^2 \mathcal{D}^\alpha u_x = 0, \quad \text{with initial condition} \quad u(x,0) =u_0(x),
    \end{equation}
    with $\alpha = 1$, and compare it with the Benjamin-Ono equation \cite{thomee1998numerical}
    \begin{equation}\label{boenn}
         u_t + uu_x - \varepsilon^2\mathcal{H}u_{xx} = 0, \qquad u(x,0) = u_0(x),
    \end{equation}
     where $\mathcal{H}$ is the Hilbert transform \cite{dutta2016convergence}. The associated initial data is given by
    \begin{equation}\label{BOsolution1}
        u_0(x) = \frac{2c\delta}{1 - \sqrt{1 - \delta^2}\cos(c\delta x)}, \qquad \delta = \frac{\pi}{cL}.
    \end{equation}
    Let $u_0$ be the initial condition associated with the fractional KdV equation \eqref{usualkdvbo}, and we use the parameters $L = 15$ and $c = 0.5$. Figure \ref{fig:timevarbo} represents that the oscillations arise after the break time $\tilde t_{c}$, where $\tilde{t}_{c} = 1/\displaystyle\max_{x\in \R}[-u_0'(x)]$. In the form of \eqref{asyptoticsol}, the asymptotic solution in the explicit form is not known in this case. However, we have observed in Figure \ref{fig:timevarbo} that the oscillations have similar asymptotic after $t>\tilde{t}_c$ as seen in \cite{miller2011zero}.

\begin{figure}
    \centering
    \includegraphics[width=0.8 \linewidth, height=7cm]{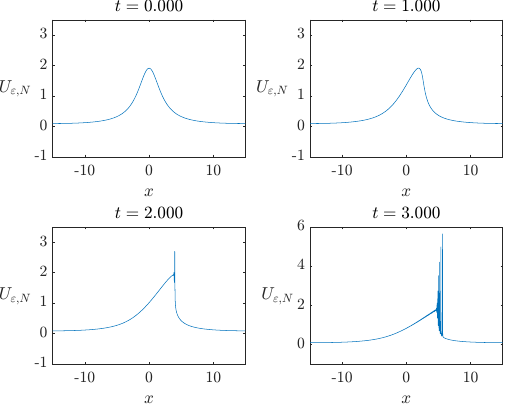}
    \caption{ The numerical approximation $U_{\varepsilon,N}$ of the fractional KdV equation \eqref{usualkdvbo} at different times $t$, and with the coefficient $\varepsilon=10^{-1}$ and for the $\alpha=1$.}
    \label{fig:timevarbo}
\end{figure}

Here are some observations on the behavior of numerical and asymptotic solutions as $\varepsilon \to 0$ for $\alpha \in [1,2]$:

\begin{itemize}
    \item Our numerical experiments reveal the presence of small oscillations even before the break time $t_c$, consistent with findings by Grava and Klein \cite{grava2007numerical}. The left plot in Figure \ref{fig:breaktime} shows oscillations near $x_c$ at $t = 0.214$, just before $t_c$, for $\varepsilon = 10^{-1}$ and $\alpha =1.999$. As $\varepsilon$ tends to 0, the oscillations diminish, and the solution converges to the solution of the limiting equation \eqref{eqn_burgers}. For $\varepsilon = 10^{-4}$, the right plot in Figure \ref{fig:breaktime} shows no oscillations before $t_c$.

    \item {Numerical simulations reveal that the oscillatory region of the numerical solution \(U_{\varepsilon,N}\) for the fractional Korteweg-de Vries (KdV) equation \eqref{fkdv} slightly extends beyond the theoretical oscillatory zone \([x^-(t), x^+(t)]\), as observed in the study by Grava and Klein \cite{grava2007numerical} for the KdV equation. Furthermore, Figure \ref{fig:errzdl} confirms that, at \(t = 0.4 > t_c\), the oscillatory region in the numerical simulations is larger than the corresponding theoretical region.} As $\varepsilon \to 0$, the numerical oscillatory zone contracts, approaching $[x^-(t), x^+(t)]$, as shown in Figure \ref{fig:epsvar}. For a fixed small $\varepsilon$, the oscillatory zone $[x^-(t), x^+(t)]$ expands with time for $t > t_c$, and Figure \ref{fig:timevar} illustrates significant oscillations for small $\varepsilon$ as time progresses past $t_c$.

    \item We have observed that decreasing of exponent $\alpha$ leads to increased in oscillations due to weaker dispersion. Figure \ref{fig:alphavar} shows that variations in $\alpha$ affect the approximate solution, also reduces the numerical oscillatory zone for fixed $\varepsilon = 10^{-1}$ and time $t = 0.5$.

\item The error estimates presented in Table \ref{tab:error_tableeps} with respect to $\varepsilon$ demonstrate that the error decreases proportionally with $\varepsilon$. This indicates that the numerical solution closely approximates the exact solution of \eqref{fkdv} even for sufficiently small values of $\varepsilon$ when $t = 0.2 < t_c$ and $t = 0.4 > t_c$. For larger values of $\varepsilon$ (of order 1), the significant dispersion eliminates the oscillatory behavior for smooth solutions, allowing the numerical solution $U_{\varepsilon,N}$ to align well with the fractional KdV equation \eqref{fkdv}.

\end{itemize}


\section{Numerical Illustrations}\label{sec5}

In this section, we aim to validate our theoretical findings, including the convergence rate for the dispersion coefficient $\varepsilon$ of order 1. We explore various examples with different values of $\alpha$ within the interval $[1,2]$ to observe the influence of $\alpha$ on the solution of the fractional KdV equation \eqref{fkdv}. The numerical solutions obtained using the spectral Galerkin scheme \eqref{eqn_fdfsgscheme} are compared with the exact solutions of the KdV equation and the Benjamin-Ono equation, corresponding to $\alpha = 2$ and $\alpha = 1$, respectively.
To ensure the robustness of our comparisons, we consider several examples from the literature \cite{dutta2016convergence,dwivedi2023stability,dutta2015convergence,thomee1998numerical}, adjusted with a fixed constant coefficient of the nonlinear term. Additionally, we verify that the integral quantities are conserved by the scheme \eqref{eqn_fdfsgscheme} as established in Lemma \ref{lem_exis}. The corresponding normalized integral quantities are defined as follows
\begin{align*}
    I^1_N := \frac{\int_{I} U_{\varepsilon,N}\,dx}{\int_{I} u_0\,dx},\quad
    I^2_N := \frac{\|U_{\varepsilon,N}\|_{L^2(I)}}{\norm{u_0}_{L^2(I)}},
    \quad I^3_N := \frac{\int_{I} \left((\mathcal{D}^{\alpha/2}U_{\varepsilon,N})^2 - \frac{1}{3}(U_{\varepsilon,N})^3\right)~dx}{\int_{I} \left((\mathcal{D}^{\alpha/2}u_0)^2 - \frac{1}{3}(u_0)^3\right)~dx}, \quad \alpha\in[1,2],
\end{align*}
{ {where $u_0$ is given periodic initial data over $I$.}
We compute the rate of convergence using the following expressions
\begin{equation}\label{rateoc}
    R =  \frac{\ln(E(N_1))-\ln(E(N_2))}{\ln(N_2)-\ln(N_1)},
\end{equation}
where $L^2$-error $E$ corresponding to scheme \eqref{eqn_fdfsgscheme} is treated as functions dependent on the number of basis trigonometric polynomials $N_1$ and $N_2$.  {We use the time step $\Delta t = 1/(N\norm{u_0}_\infty)$ in the subsequent numerical simulations.}
We begin with the following example.

\subsection*{Example 5.1}
    We consider the classical KdV equation \cite{dwivedi2023stability,dutta2015convergence} 
    \begin{equation}\label{kdveqn1}
        u_t+uu_x + u_{xxx} = 0, \qquad u(x,0) = u_0(x),
    \end{equation}
    which is close to \eqref{fkdv} with $\alpha = 1.999$ up to constant coefficients. The one-soliton solution \cite{dutta2015convergence,dwivedi2023stability} of the KdV equation \eqref{kdveqn1} is given by
    \begin{equation}\label{onesol}
        u(x,t) = 9\left(1 - \tanh^2\left(\sqrt{\frac{3}{2}}(x - 3t)\right)\right).
    \end{equation}
    Let $u_0(x) = u(x,0)$ be the initial condition associated with the fractional KdV equation \eqref{usualkdvbo} with $\alpha=1.999$ and the KdV equation \eqref{kdveqn1}. We compute the solution of \eqref{fkdv} at time $t = 2$ using the fully discrete scheme \eqref{eqn_fdfsgscheme}, with the dispersion coefficient $\varepsilon = 1$, and compare it with the solution \eqref{onesol} of the KdV equation \eqref{kdveqn1}. 
 
\begin{table}
  \centering
  \begin{tabular}{||c|c|c|c|c|c||}
        \hline
 N & $E$ & $R$ & $I^1_{N}$ & $I^2_{N}$ & $I^3_{N}$ \\
 \hline
 \hline
   128  & 8.15e-04 & &  1.00 & 1.00 & 1.00  \\
  &  & 2.53  & & &\\
  256 & 1.40e-04 &&1.00 & 1.00 & 1.00 \\
  &  &  2.00 & &&\\
 512 &  3.49e-05& &1.00 & 1.00 & 1.00  \\
  &  & 2.00  &&& \\
 1024 & 8.72e-06 & &1.00 & 1.00 & 1.00 \\
   &  & 1.98&&& \\
 2048 & 2.20e-06 &&1.00 & 1.00 & 1.00 \\
\hline
    \end{tabular}
   \caption{Errors $E$, rate of convergence $R$ and integral quantities $I^i_N$, $i=1,2,3$  with $\alpha= 1.999$, $T=2$ and initial data $u_0 = u(x,0)$ given by \eqref{onesol}.}
   \label{tab:TableKdVCN}
\end{table}
\begin{figure}
    \centering
    \includegraphics[width=0.8 \linewidth, height=8cm]{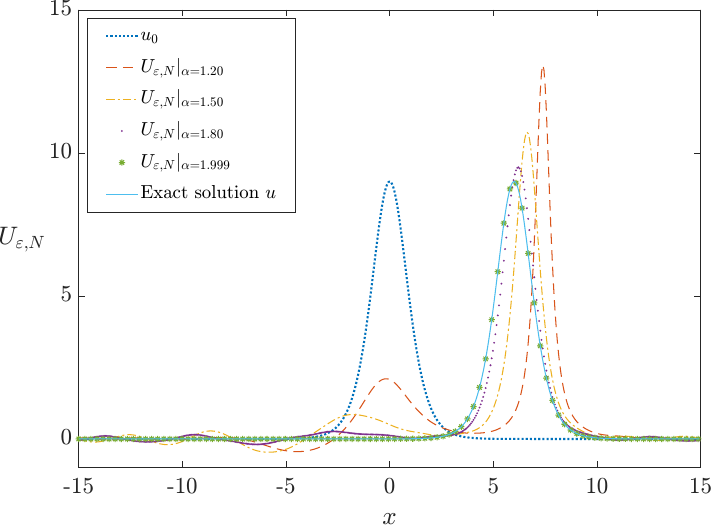}
    \caption{Numerical approximation $U_{\varepsilon,N}$ with $N=512$ of the fractional KdV equation \eqref{fkdv} for various values of $\alpha$ at time $t = 2$, compared with the asymptotic solution $u$ of the KdV equation \eqref{kdvellson}.}
    \label{fig:kdv2onsol}
\end{figure}
Figure \ref{fig:kdv2onsol} shows the numerical approximation $U_{\varepsilon,N}$ of the fractional KdV equation \eqref{fkdv} for various values of $\alpha$ at time $t = 2$. It compares the numerical results with the asymptotic solution $u$ of the KdV equation \eqref{kdvellson}, highlighting the accuracy of the spectral Galerkin scheme \eqref{eqn_fdfsgscheme} as $\alpha$ approaches $2$. The results suggest how varying $\alpha$ affects the dispersive properties of the solutions. The figure illustrates that as $\alpha$ gets closer to $2$, the numerical solution converges to the usual KdV solution \eqref{onesol}, showing a convergence pattern and confirmed by Table \ref{tab:TableKdVCN}. Moreover, Table \ref{tab:TableKdVCN} depicts that the integral quantities $I^i_N,~i=1,2,3$ are conserved in the discrete setup.
\subsection*{Example 5.2}
    We consider the fractional KdV equation \eqref{usualkdvbo} with $\alpha = 1.01$ and $\varepsilon = 1$, and compare it with the Benjamin-Ono equation \cite{thomee1998numerical} given by \eqref{boenn}.
    The one-soliton solution of the Benjamin-Ono equation is given by \cite{dutta2016convergence,thomee1998numerical}
    \begin{equation}\label{BOsolution}
        u_1(x,t) = \frac{2c\delta}{1 - \sqrt{1 - \delta^2}\cos(c\delta(x - ct))}, \qquad \delta = \frac{\pi}{cL}.
    \end{equation}
    Let $u_0 = u_1(x,0)$ be the initial condition associated with the fractional KdV equation \eqref{usualkdvbo}, using the parameters $L = 15$ and $c = 0.25$.

\begin{table}
  \centering
  \begin{tabular}{||c|c|c|c|c|c||}
        \hline
 N & $E$ & $R$ & $I^1_{N}$ & $I^2_{N}$ & $I^3_{N}$ \\
 \hline
 \hline
   128  & 9.36e-09 & &  1.00 & 1.00 & 1.00  \\
  &  & 2.00  & & &\\
  256 & 2.34e-09 &&1.00 & 1.00 & 1.00 \\
  &  &  2.00 & &&\\
 512 & 5.84e-10 & &1.00 & 1.00 & 1.00  \\
  &  & 1.99  &&& \\
 1024 & 1.46e-10 & &1.00 & 1.00 & 1.00 \\
   &  & 1.99&&& \\
2048 & 3.67e-11 &&1.00 & 1.00 & 1.00 \\
\hline
    \end{tabular}
   \caption{Errors $E$, rate of convergence $R$ and integral quantities $I^i_N$, $i=1,2,3$  with $\alpha= 1.01$, $T=20$ and initial data $u_0 = u_1(x,0)$ given by \eqref{BOsolution}.}
   \label{tab:TableboCN}
\end{table}

Figure \ref{fig:bo1onsol} displays the numerical approximation $U_{\varepsilon,N}$ of the fractional KdV equation \eqref{fkdv} for $\alpha\approx1$ at time $t=20$, alongside the classical solution \eqref{BOsolution} of the Benjamin-Ono equation \eqref{boenn}. The second order convergence rate in time is obtained and presented in Table \ref{tab:TableboCN}. 
\begin{figure}
    \centering
    \includegraphics[width=0.8 \linewidth, height=8cm]{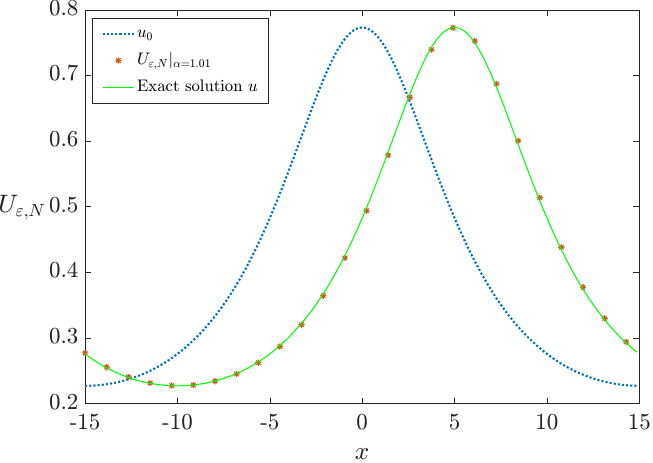}
    \caption{Numerical approximation $U_{\varepsilon,N}$ with $N=512$ of the fractional KdV equation \eqref{fkdv} for $\alpha=1.01$ at time $t = 20$, compared with the asymptotic solution $u$ of the Benjamin-Ono equation \eqref{BOsolution}.}
    \label{fig:bo1onsol}
\end{figure}

\subsection*{Example 5.3}
    We consider the fractional KdV equation \eqref{usualkdvbo} with $\alpha = 1.5$ and $\varepsilon = 1$. The initial condition is set as $u_0(x) = 0.5\sin(x)$ for $x \in [-\pi, \pi]$. The approximate solution is obtained using the numerical scheme \eqref{eqn_fdfsgscheme} and compared with a reference solution computed using a higher grid with $N = 2^{16}$ at $t=2$. Table \ref{tab:Tablefrac1_5} shows the convergence rates, confirming that the method achieves second order accuracy in time, validating the theoretical convergence results.

\begin{table}
  \centering
  \begin{tabular}{||c|c|c|c|c|c||}
        \hline
 N & $E$ & $R$ & $I^1_{N}$ & $I^2_{N}$ & $I^3_{N}$ \\
 \hline
 \hline
   128  & 2.64e-06 & &  1.01 & 1.00 & 1.00  \\
  &  & 1.99  & & &\\
  256 & 6.64e-07 &&1.00 & 1.00 & 1.00 \\
  &  &  2.00 & &&\\
 512 & 1.66e-07 & &1.00 & 1.00 & 1.00  \\
  &  & 1.00  &&& \\
 1024 & 4.14e-08 & &1.00 & 1.00 & 1.00 \\
   &  & 2.01&&& \\
2048 & 1.02e-08 &&1.00 & 1.00 & 1.00 \\
\hline
    \end{tabular}
   \caption{Errors $E$, rate of convergence $R$ and integral quantities $I^i_N$, $i=1,2,3$  with $\alpha= 1.5$, $T=2$ and initial data $u_0 = 0.5\sin x$.}
   \label{tab:Tablefrac1_5}
\end{table}

\section{Concluding Remarks}\label{sec6}
In this study, we have developed a structure-preserving Fourier spectral Galerkin (FSG) scheme which conserves mass, momentum, and energy, and demonstrated that the fully discrete Crank-Nicolson (CN) FSG scheme is stable and convergent. Additionally, we provided a constructive proof for the existence and uniqueness of the solution to the fractional KdV equation using compactness arguments. The proposed scheme achieves optimal spectral accuracy for periodic initial data in $H^r$, and exponential accuracy for analytic initial data, resulting in significantly reduced computational time during numerical simulations compared to our previous works.

Moreover, we investigated the zero-dispersion limit of the fractional KdV equation. Through extensive numerical simulations, we confirmed that the proposed CN-FSG scheme accurately captures the oscillations generated by small dispersion. We also examined the behavior of the numerical solutions as $\alpha$ varies within $[1,2]$, and evaluated the performance of the scheme for small values of $\varepsilon$ beyond the critical time \(t_c\).
However, further work remains, both theoretically and numerically, to better understand the asymptotic behavior of the zero dispersion limit of the fractional KdV equation beyond \(t_c\).

\section*{Acknowledgements and conflicts of interest} 
The authors declare that they have no known competing financial interests that could have appeared to influence the work reported in this paper. Furthermore, no data was used for the research described in the article.
\section*{Declarations}

\subsection*{Funding}
This research received no external funding.

\subsection*{Data Availability}
Data sharing is not applicable to this article as no datasets were generated or analyzed during the current study.

\subsection*{Competing Interests}
The authors declare that they have no competing interests.

\subsection*{Author Contributions}
Mukul Dwivedi: Conceptualization, methodology, visualization, writing original draft, numerical experiments.\\
Tanmay Sarkar: Supervision, conceptualization, validation, writing, review, and editing.

\begin{appendices}
\section{Appendix} 
\renewcommand{\thelemma}{A.\arabic{lemma}} 
 We collect a number of elementary technical results that were used earlier in the article.
\begin{lemma}[Inverse and Sobolev-type inequalities for trigonometric polynomials]\label{A.1}
Let $f\in V_N$. Then:
\begin{enumerate}
    \item (Inverse inequality) For the $L^2$-norm, we have
    \begin{equation}
        \|f_x\| \leq N\|f\|.
    \end{equation}
    \item (Sobolev inequality) We have
    \begin{equation}
        \|f\|_\infty \leq \sqrt{\frac{10}{3}}\|f\|_1.
    \end{equation}
    \item (Sobolev embedding) For any $r \geq 0$, we have
    \begin{equation}
        \|f\| \leq \|f\|_r.
    \end{equation}
\end{enumerate}
\end{lemma}

\begin{proof}  {Note that for $f\in V_N$ we have
\begin{equation*}
 f(x) := \sum_{k=-N}^{N} \hat{f}(k) e^{ikx}~ \text{ and }~ f_x(x) := \sum_{k=-N}^{N} ik\hat{f}(k) e^{ikx}.
\end{equation*}}
\begin{enumerate}
    \item  {For the inverse inequality, we compute using the Parseval's identity
    \begin{align*}
        \|f_x\|^2 &= \sum_{k=-N}^N |k|^2 |\hat{f}(k)|^2 \leq N^2 \sum_{k=-N}^N |\hat{f}(k)|^2 = N^2\|f\|^2.
    \end{align*}
    Taking square roots gives the result.}

    \item  {For the Sobolev inequality, we first bound the sup-norm and use the Cauchy-Schwarz inequality
    \begin{align*}
        \|f\|_\infty &\leq \sum_{k=-N}^N |\hat{f}(k)| 
        = \sum_{k=-N}^N \frac{(1 + k^2)^{1/2}|\hat{f}(k)|}{(1 + k^2)^{1/2}} \\
        &\leq \left(\sum_{k=-N}^N \frac{1}{1 + k^2}\right)^{1/2} \left(\sum_{k=-N}^N (1 + k^2) |\hat{f}(k)|^2\right)^{1/2}.
    \end{align*}
    The second factor is exactly $\|f\|_1$. We now bound the summation factor
    \begin{align*}
        S_N &:= \sum_{k=-N}^N \frac{1}{1 + k^2} = 1 + 2\sum_{k=1}^N \frac{1}{1 + k^2} \leq 1 + 2\sum_{k=1}^\infty \frac{1}{1 + k^2} = 1+2\left(\frac{\pi\coth \pi -1}{2}\right)\\
        &= \pi \coth \pi < \frac{10}{3}.
    \end{align*}
 Hence, we have 
 $$\|f\|_\infty \leq \sqrt{\frac{10}{3}}\|f\|_1.$$}
    \item  {For the Sobolev embedding ($r \geq 0$)
    \begin{align*}
        \|f\|^2 &= \sum_{k=-N}^N |\hat{f}(k)|^2 = \sum_{k=-N}^N (1 + k^2)^0 |\hat{f}(k)|^2 \leq \sum_{k=-N}^N (1 + k^2)^r |\hat{f}(k)|^2 = \|f\|_r^2,
    \end{align*}
    since $(1 + k^2)^r \geq 1$ for all $k$ when $r \geq 0$.}
\end{enumerate}
\end{proof}






\end{appendices}
 

\end{document}